\documentclass[10pt, runningheads]{llncs}
\usepackage{amsmath}
\usepackage{amssymb}
\usepackage{color}
\usepackage{fancyhdr}
\usepackage{graphicx}
\usepackage{hyperref}
\usepackage{xcolor}
\usepackage[polutonikogreek, english]{babel}

\newcommand{\pvn}{\par\vspace{1ex}\noindent}

\makeatletter
\renewcommand\subsubsection{\@startsection{subsubsection}{3}{\z@}%
                       {-18\p@ \@plus -4\p@ \@minus -4\p@}%
                       {0.5em \@plus 0.22em \@minus 0.1em}%
                       {\normalfont\normalsize\bfseries\boldmath}}
\makeatother
\begin{document}

\title{Phenomenology of diagrams\\ in Book II of the \emph{Elements}}

%
%\titlerunning{Abbreviated paper title}
% If the paper title is too long for the running head, you can set
% an abbreviated paper title here
%
\author{Piotr Błaszczyk}%\inst{1}\orcidID{0000-0002-3501-3480} %\and \\ Anna Petiurenko\inst{2}\orcidID{0000-0002-0196-6275} %\and
%Third Author\inst{3}\orcidID{2222--3333-4444-5555}

%
\authorrunning{Piotr Błaszczyk} %Anna Petiurenko}
% First names are abbreviated in the running head.
% If there are more than two authors, 'et al.' is used.
%
\institute{Institute of Mathematics, Pedagogical University of Cracow, Podchorazych 2, Cracow, Poland\\
\email{piotr.blaszczyk@up.krakow.pl}} 
%\and Institute of Mathematics, Pedagogical University of Cracow, Podchorazych 2, Cracow, Poland\\
%\email{piotr.blaszczyk@up.krakow.pl} %;\ anna.petiurenko@up.krakow.pl}
%\url{https://matematyka.up.krakow.pl/}
%\and
%ABC Institute, Rupert-Karls-University Heidelberg, Heidelberg, Germany\\
%\email{\{abc,lncs\}@uni-heidelberg.de}

%
%\maketitle              % typeset the header of the contribution
%
{\def\addcontentsline#1#2#3{}\maketitle}
\begin{abstract} 
In this paper, we provide an interpretation of Book II  of the \textit{Elements} from the perspective of 
figures which are represented and not represented on the diagrams. We show that Euclid's reliance on figures not represented on the diagram  is a proof technique which enables to turn his diagrams II.11--14 into ideograms of a kind.

  We also discuss interpretations of Book II developed by J. Baldwin and A. Mueller, L. Corry, D. Fowler,
R. Hartshorne,  I. Mueller, K. Saito, and the so-called geometric algebraic interpretation in B. van der Waerden's version.

\end{abstract}
%
%

%\section{Euclid's unrepresented figures}
\tableofcontents
\newpage\section{Introduction}
Proposition II.1 of  Euclid's \emph{Elements} states that ``the rectangle contained by A, BC is equal to the rectangle contained by A,\,BD, by A,\,DE, and, finally, by A,\,EC", given BC is cut at  D and E.\footnote{All English translations of the \textit{Elements} after (Fitzpatrick 2007). Sometimes we slightly modify
Fitzpatrick's version by skipping interpolations, most importantly, the words related
to addition or sum. Still, these amendments are easy to verify, as this edition is
available on the Internet, and also provides the Greek text and diagrams of the classic Heiberg
edition.} In algebraic stylization, it is conveyed by the formula $A\cdot BC= A\cdot BD+A\cdot DE+A\cdot EC$, where $BC=BD+DE+EC$.
In  modern theory of rings, this formula  %$A\cdot (BD+DE+EC)= A\cdot BD+A\cdot DE+A\cdot EC$ 
is simplified to $a(b+c)=ab+ac$ and represents  distributivity  of multiplication over addition. 
 In algebra, however, it is an axiom, therefore, it seems unlikely that Euclid managed to prove 
it, even  in a geometric disguise. Moreover, if we apply algebraic formulea to read proposition II.1, it appears that
the equality $a(b+c+c)=ab+ac+dc$ is both the starting-point and the conclusion of the proof. Yet, there is some  in-between in Euclid's proof. What is this residuum  about? Although an algebraic formalization easily interprets the thesis of the proposition, it does not help to reconstruct its proof.

The above interpretation is an extreme. Usually, instead of the multiplication $A\cdot BC$ the term $A.BC$ is considered, which stands for a rectangle with  sides $A$ and $BC$. Then, it is assumed that the equality $A.BC=A.BD+A.DE+A.EC$  is approved by the accompanying diagram. In this way, a mystified role of Euclid's diagrams substitute detailed analyses of his proofs.

David Fowler, for example, ascribes to Euclid's diagrams in Book II not only the  power  of proof makers. In his view, they  summarize both subjects and proofs of propositions:  ``The subject of each proposition is best conveyed by its figure (and it is these figures, not what is made of them in their enunciations or proofs, that will enter my proposed reconstruction)" (Fowler 2003, 66).

Accordingly, he presents propositions II.1--8 as arguments based on diagrams: ``In all of these propositions, almost all of the text  of the demonstrations concerns the construction of the figures, while the substantive content of each enunciation is merely  read off from constructed figures, at the end of the proof (Fowler 2003, 69).
In propositions II.9--10, Euclid studies  the use of the Pythagorean theorem, and the accompanying diagrams represent right-angle triangles rather then squares descried on their sides. Yet, to buttress his interpretation,  Fowler provides  alternative proofs, as he believes Euclid basically  applies  ``the technique of dissecting squares". In his view,  Euclid's proof technique is very simple:
``With the exception of implied uses of I47 and 45, Book II is virtually self-contained in the sense  that it only uses straightforward manipulations of lines and squares of the kind assumed without comment by Socrates in the \textit{Meno}"(Fowler 2003,  70). 

Fowler is so focused on dissection proofs that he cannot spot what actually is and what is not depicted on Euclid's diagrams. As for proposition II.1,  there is clearly no rectangle contained by \textit{A} and \textit{BC}, although there is a rectangle with vertexes \textit{B,\,C,\,H,\,G} (see Fig. \ref{figII1}).  Indeed, all throughout Book II Euclid deals with figures which are not represented on diagrams. Finally, in propositions II.11, 14 they appear to be a tool to establish results concerning figures represented on the diagrams. 

The plan of the paper is as follows. Section \S\,2 provides an overview of Book II, specifically we discuss four groups of propositions, II.1--3, II.4--8, II.9--10, and II.11--14, in terms of applied  proof techniques.
In section \S\,3, we analyze basic components of Euclid's propositions: lettered diagrams, word patterns, and the concept of \textit{parallelogram contained by}. In section \S\,4, we scrutinize propositions II.1--4 and introduce symbolic schemes of Euclid's proofs. In section \S\,5, we introduce non-geometrical rules which aim to explain relations between figures which are represented (visible) and not represented (invisible) on Euclid's diagrams. In section \S\,6, we analyze the use of propositions II.5--6 in II.11,\,14 to demonstrate how the technique of invisible figures enables to establish relations between visible figures. 

Throughout the paper we confront specific aspects of our interpretation with
readings of Book II by scholars such as  J. Baldwin and A. Mueller, L. Corry, D. Fowler,
R. Hartshorne,  I. Mueller,  and B. van der Waerden. In section \S\,7, we discus overall interpretations of Book II.
Finally, in section \S\,8, we discuss proposition II.1 from the perspective of Descartes's lettered diagrams. 
We show that there is a germ of algebraic style in Book II, however, it has not been developed further in the \textit{Elements}.

\section{Overview of Book II}

\subsection{Three groups of propositions}
 Book II of the \textit{Elements} consists of two definitions and fourteen propositions. The first definition introduces the term \textit{parallelogram contained by}, the second -- \textit{gnomon}.
  All parallelograms considered  are rectangles and squares, and indeed there are two basic concepts applied throughout Book II, namely,  \textit{rectangle contained by}, and \textit{square on}, while the   \textit{gnomon} is used only in propositions II.5--8. 
  
Considering the results, proof  techniques, and word  and diagrammatic patterns, we distinguish three groups of propositions: 1--8, 9--10, 11--14.

 II.1--8 are lemmas. II.1--3 introduce a specific  use of the terms \textit{squares on} and \textit{rectangles contained by}.  II.4--8 determine the relations between squares and rectangles resulting from dissections of bigger squares or rectangles. Gnomons play a crucial role in these results. Yet, from II.9 on, they are of no use.

 Propositions II.9--10 apply the Pythagorean theorem for combining squares. To this end, Euclid considers right-angle triangles sharing a hypotenuse and equates squares built on their legs. Although these results could be obtained by dissections and the use of gnomons, proofs based on I.47 provide new insights.
  
   II.11--14  present substantial geometric results. Their proofs are based on the  lemmas II.4--7, and the use of the Pythagorean theorem in the way introduced in II.9--10.

   Propositions II.1--3, when viewed from the perspective of deductive structure, seem redundant. Yet, we consider them as introducing the  basics of the technique developed  further in II.4--8. Similarly, II.9--10  also seem redundant -- when viewed from this perspective. Yet, they introduce a technique of applying the Pythagorean theorem.

  \subsubsection{Ian Mueller on the structure of Book II}  In regard to the structure of Book II, Ian Mueller writes: ``What unites all of book II is the methods employed: the addition and subtraction of rectangles and squares to prove equalities  and the construction of rectilinear areas satisfying given conditions. 1--3 and 8--10 are also applications of these methods; but why Euclid should choose to prove exactly those propositions does not seem to be fully explicable" (Mueller, 2006, 302). 
  
  Our comment on this  remark is simple: the perspective of deductive structure, elevated by Mueller to the title of his book, does not cover propositions dealing with  technique. Mueller's perspective, as well as his Hilbert-style reading of the \textit{Elements}, results in a distorted, though comprehensive overview of the \textit{Elements}. Too many propositions  do not find their place in this deductive structure of \textit{the Elements}. While interpreting  the \textit{Elements}, Hilbert applies his own techniques, and, as a result, skips the propositions which specifically develop Euclid's technique, including the use of the compass. Furthermore,   in the \textit{Grundalgen}, Hilbert does not provide any proof of the Pythagorean theorem, while in our interpretation it is both a crucial result (of Book I) and  a  proof technique (in Book II).\footnote{The Pythagorean theorem plays a role in Hilbert's models, that is, in his meta-geometry.}
   
    In  modern mathematics, there are many important results  concerning proof technique. The transfer principle relating  standard and non-standard analysis, is a model example. Hilbert's proposition  that the equality of polygons built on the concept of dissection and Euclid's theory of equal figures do not produce equivalent results could be another example. Viewed from that perspective, II.9--10 show how to apply I.47 instead of gnomons to acquire the same results. II.1--3 introduce a specific use of the terms \textit{square on} and  \textit{rectangle contained by} which Mueller ignores in his analysis of Book II (see \S\,7 below).

\subsection{Basic geometric results}

 In regard to geometric results,  II.11 provides  the so-called golden ratio construction. It is a crucial step in  the cosmological plan of the \textit{Elements}, namely -- the  construction of a regular pentagon and finally, a dodecahedron. The respective justification builds on II.6.

  II.12,\,13 are what we  recognize as the cosine rule for  the obtuse and acute triangle respectively. The former proof begins with a reference to II.4, the later -- with a reference to II.7.

\begin{figure}
\begin{center}
\includegraphics[scale=1]{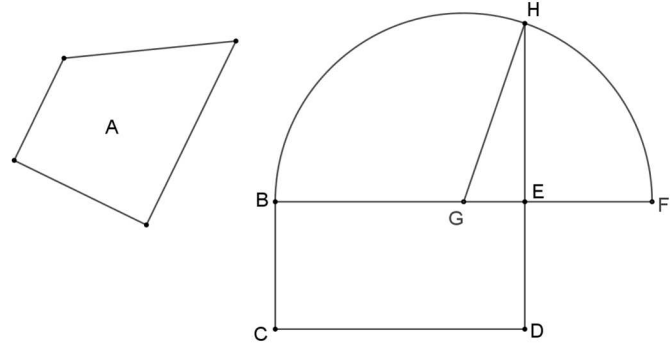}
\caption{\textit{Elements}, II.14.} \label{figII14}
\end{center}
\end{figure}  
  
 In II.14, Euclid shows how to square a polygon.
 This construction  crowns the theory of equal figures developed in propositions I.35--45; see (Błaszczyk 2018).  In Book I, it involved  showing how  to build a parallelogram equal to a given polygon. In II.14, it is already assumed that the reader knows how to transform a  polygon into an equal rectangle. The justification of the squaring of a polygon  begins with a reference to II.5.

 As for the proof technique, in II.11--14, Euclid combines the results of II.4--7 with the Pythagorean theorem by adding or subtracting squares described on the sides of  right-angle triangles. 
Thus, II.5, 6 share the same scheme: 
 \[rectangle + \textit{square on}\ A= \textit{square on}\ B.\]  
 
 When applied, 
a right-angle triangle with a hypotenuse  \textit{B} and legs  \textit{A}, \textit{C} is considered. Then, by  I.47, 
\[rectangle + \textit{square on}\ A= \textit{square on}\ A + \textit{square on}\ C.\]
 
 By subtraction from  both sides of the \textit{square on}\ A, the equality characterizing II.11 and II.14 is obtained
  \[rectangle= \textit{square on}\ C.\]
  
  On the other hand, II.4,\,7 share another scheme. In II.4, it is as follows
  \[\textit{square on}\ A= \textit{square on}\ B+ \textit{square on}\ C+ 2 rectangles. \]
  
  Addition to both sides another square gives the equality
   \[\textit{square on}\ A+ \textit{square on}\ D = \textit{square on}\ B+ \textit{square on}\ D+ \textit{square on}\ C+ 2 rectangles. \]

 \textit{D}  is a leg of two right-angle triangles: the first   with another leg \textit{A} and a hypotenuse \textit{E}, the second with   leg \textit{B} and hypotenuse  \textit{F}. Thus,  equality obtains\footnote{The term \textit{2rectangles} can be interpreted as $-2FC\cos(\textit{angle between sides F and C})$. See section \S\,6.2 below.}
 \[\textit{square on}\ E= \textit{square on}\ F + \textit{square on}\ C + 2 rectangles. \]
 
The use of II.7 starts with the equation
\[\textit{square on}\ A+ \textit{square on}\ B= \textit{square on}\ C+ 2rectangles.\]

\begin{figure}
\begin{center}
\includegraphics[scale=1]{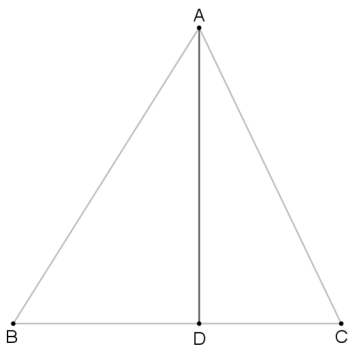}
\caption{\textit{Elements}, II.13.} \label{figII13}
\end{center}
\end{figure}

A square on the height \textit{AD }of the triangle \textit{ABC} is added to its both sides; see Fig. \ref{figII13}.  The rest of the proof  II.13 proceeds in a similar way.\footnote{Using the names of lines in the diagram, $BC^2+DC^2=BD^2+2BC\times DC$. By adding $AD^2$ to the both sides of this equation, we obtain $BC^2+AC^2=BA^2+2BC\times DC$.}

  Significantly, the accompanying diagram, next to sides of the triangle, features only the height \textit{AD}. Thus, the point \textit{D }represents  the way the side \textit{BC} is cut, namely \textit{at random}. In this way, it makes a reference to II.7. The line \textit{AD} represents the side of the square that is added to  both sides of the above equation.  
  All that illustrates our thesis,  that although in II.11--14 Euclid relates squares and rectangles, the accompanying diagrams depict only the respective sides of the figures involved. 
% Only by this specific use of the Pythagorean theorem,\textit{rectangles contained by} are made equal to simple %figures represented on the diagrams.
  
%\end{document}
  \subsection{In the search  of  patterns}
 References to II.4--7 follow a similar word pattern: \textit{For since the straight-line ... has been cut ...}\,. Then comes a specification   how  new points are placed on the straight-line:  \textit{in half} (\textit{equally}),  \textit{at random} (\textit{unequally}), or a new line is added to the original one. Finally, a relation between resulting squares and rectangles established in the refereed proposition is reformulated due to new names of points.

 For example, the respective part of II.14 is this: ``since the straight-line BF has been cut
equally at G, and unequally at E, the rectangle contained by BE, EF, together with the square on EG, 
is thus  (\foreignlanguage{polutonikogreek}{>'ara})
equal to the square on GF".
 
While II.5 states: 
 ``For let any straight-line AB have been cut equally at C and unequally at D. I say that the rectangle contained
by AD, DB, together with the square on CD, is equal to the
square on CB".

Clearly, these statements differ only in names of points (for comparison, see  Fig. \ref{figII14} and Fig. \ref{figII5}).

In propositions II.1--8, depending on a distribution of cut points, a variety of squares and rectangles appears on the accompanying diagrams. Next to those represented on the diagrams, Euclid refers also to figures  which are not represented on the diagrams; let us call them  invisible figures. These are  \textit{squares on}  and \textit{rectangles contained by} some  lines. Although the respective lines are represented on the diagrams, the related squares and rectangles are not. Indeed, whereas these invisible figures occur in the statements of propositions, Euclid's proofs usually start with figures which are represented on the diagrams.

We recognized the following pattern in  the procedures Euclid adopts in propositions II.1--8. Firstly, contrary to Book I, the \emph{diorisomos} part of the proposition refers to figures not represented on the diagram. Secondly, in the \emph{kataskeu\={e}},  a geometric machinery is applied to construct figures represented on the diagram. Thirdly, in the \emph{apodeixis}, a  relation between  figures which are represented and not represented on the diagram is determined.  It is achieved by visual evidence,  substitution rules, and the renaming of figures. In sections \S\,4--5, we will detail these issues, while section \S\,3 is dedicated to Euclid's use of the terms \textit{square on} and  \textit{rectangle contained by}.

\begin{figure}
\begin{center}
\includegraphics[scale=1]{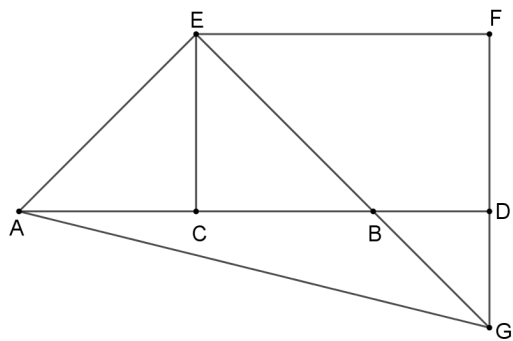}
\caption{\textit{Elements}, II.10.} \label{figII10}
\end{center}
\end{figure}

 As we proceed from II.1 to II.8, Euclid's diagrams get more complicated: they depict more and more  squares and rectangles. Then, in propositions II.9--10, they gain a new clarity. 
 Indeed,  II.9--10 explore the Pythagorean theorem in equating groups of squares, yet, the accompanying  diagrams do not depict these squares. For example,  
II.10 reads: ``For let any straight-line AB have been cut in half at C, and let any straight-line BD have been added
to it straight-on. I say that the squares on  AD, DB is double the squares on AC, CD" (see Fig. \ref{figII10}). 
 
 Since  the right-angle triangles \textit{ADG} and \textit{AEG} share the hypotenuse \textit{AG}, {squares on} the legs of the first triangle,  $AE$, $EG$, and on the legs of the second,  $AD, DG$,  are equal.\footnote{In II.9, Euclid does the same trick.} The squares on  
 $AE$ are $AC^2, CE^2$; the squares on $EG$ are $FE^2, FG^2$; and the squares on the legs of the second triangle are 
 $AD^2, DG^2$.\footnote{The term $AC^2$ stands for the phrase \textit{the square on the straight-line AC}.} Finally, Euclid determines equalities between squares to  get the \textit{doubles} of the squares, such as $2FE^2$ which is to be equal to the  square on \textit{EG}. Although the equality of line segments, $FG=EF$, is derived from I.6,  the equality of respective squares, $FG^2=FE^2$, is based on an implicit rule: ``since FG is equal to EF, the one on FG is equal to the one on EF".   In our interpretation,  it is one of substitution rules  discussed  in section \S\,5. 
 
  Notice, that to justify the  equality $FG^2=FE^2$, Euclid does not refer to the diagram. In fact, the accompanying diagram does not  depict any square.  It is also not the case, that this conclusion is taken for granted, since he provides a reason. Indeed, throughout Book II, Euclid reiterates the argument: \textit{since X is equal to Y, the square on X is equal to the square on Y}.

 Interestingly, the results II.9--10 could be obtained by dissection of the square on \textit{AB}, in the case of II.9,
 and the square on \textit{AD}, in the case of II.10, and then by the use of gnomons in a way similar to the proofs of II.5--8.   Therefore, we  view II.9--10 as introducing a new technique which combined with the results of lemmas II.4--7 is used in the proofs II.11--14.

 The word pattern of references we identified above finds its diagrammatic counterparts in II.11--14.
Starting from II.11, \textit{squares on} are represented by  lines, \textit{rectangles contained by} -- by lines with cut points. %In diagrams accompanying propositions II.11--14, references to II.4--7 are represented by a straight-line and respective cut points.
In  II.11--14,  diagrams  look like ideograms rather than 
a simple composition of lines introduced throughout construction steps.

 Let us focus on II.14 (see Fig. \ref{figII14}). 
The figure \textit{A} and the rectangle \textit{BCDE} refer to the theory of equal figures developed in Book I. 
The exposition of the  line \textit{BF} and the cut points \textit{G}, \textit{E}, refer to II.5. The semicircle \textit{BHF} and the radius \textit{GH} represent a new construction. The line \textit{EH} represents the conclusion of the squaring of a polygon process.  

\begin{figure}
\begin{center}
\includegraphics{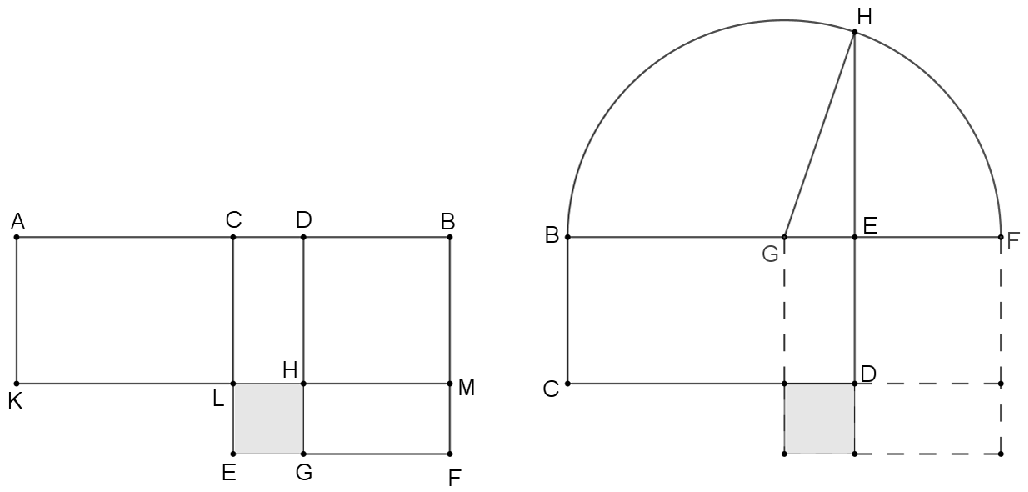}
\caption{The scheme of application of II.5 in II.14.} \label{figII5_14}
\end{center}
\end{figure}

In Fig. \ref{figII5_14}, we represent the way how II.5 is applied in II.14. 
The scheme of this application is as follows:\footnote{In section \S\,6.2, we present the scheme of the actual Euclid's proof, which applies a few substitutions such as: $GF=GH$, then $GH^2=GF^2$.} Since the line \textit{BF} is cut in half at \textit{G}, and at random at \textit{E}, by II.5, the  equality obtains 
\[rectangle\ BD + \textit{square on}\ GE =  \textit{square on}\ GF.\]

Then, by I.47: 
\[\textit{square on}\ GF=\textit{square on}\ GE + \textit{square on}\ HE.\]

Thus 
\[rectangle\ BD + \textit{square on}\ GE = \textit{square on}\ GE + \textit{square on}\ HE .\]

By subtraction, the final equality obtains

\[rectangle\ BD  =  \textit{square on}\ HE .\]

Yet, instead of all these auxiliary lines that evoke the   II.5 construction, Euclid's original diagram is much simpler: it refers to II.5 by the way points \textit{G} and \textit{E} are located on the line \textit{EF}.
Note however, that Euclid reached this clarity  due to the specific use of the term \textit{square on}.
Instead of the square on \textit{GH}, he considers square on \textit{GF}, instead of the square on \textit{GE}, he considers another square -- the gray one in Fig. \ref{figII5_14}, on the right.  And that is why, in II.5, although he mentions the square on \textit{CD}, he considers the
square \textit{LEGH} -- the gray one in Fig.  \ref{figII5_14}, on the left.

In sum, from the perspective of diagrams,  Book II applies figures which are represented and not represented on the diagrams. These are \textit{squares on} and \textit{rectangles contained by}. II.9--10  apply line segments instead of \textit{squares on}.  II.11--14, besides lines representing \textit{squares on}, apply also line segments with cut points   instead of \textit{rectangles contained by}.    

\section{Building-blocks of Euclid's propositions of Book II}

  \subsection{Individual vs  abstract components of Euclid's diagrams}

There are two components of Euclid's proposition: the text  and the lettered diagram. The Greek text is linearly  ordered  -- sentence follows sentence, from left to  right,  and from top to bottom. Diagrams consist of  
straight-lines  and circles. The capital letters on the diagrams are located next to points; they name the ends of line segments, intersections of lines, or  random points.

  The text of  the proposition is a schematic composition made up of  six parts: \emph{protasis} (stating the relations among geometrical objects by means of abstract and technical terms), \emph{ekthesis} (identifying objects of \textit{protasis} with lettered objects), \emph{diorisomos} (reformulating \textit{protasis} in terms of lettered objects), \emph{kataskeu\={e}}   (a construction part which introduces auxiliary lines exploited in the proof that follows),  \emph{apodeixis} (proof, which usually proves  the \emph{diorisomos}' claim),  \emph{sumperasma} (reiterating \textit{diorisomos}).  References to axioms, definitions, and previous propositions  are made via  the technical terms and phrases applied in \emph{prostasis}. 
  
  \begin{figure}
\begin{center}
\includegraphics[scale=1.2]{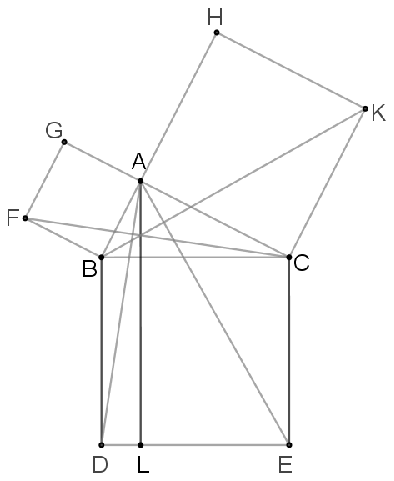}
\caption{\textit{Elements}, I.47.} \label{figI47}
\end{center}
\end{figure}

  In fact, it is the received account of Euclid's propositions. For example, in proposition I.47, the \emph{protasis}, \emph{ekthesis}, and \emph{diorisomos} 
  are as follows (numerals in square brackets added, here they stand for subsequent  parts of the proposition):
  
  [1] ``In a right-angle triangle, the square on the side
subtending the right-angle is equal to squares on the sides surrounding the right-angle.

[2] Let ABC  be a right-angled triangle having the right-angle BAC. 

[3] I say that the square on BC is equal to the
squares on BA and AC."

Indeed, the triangle \textit{ABC} as well as the squares constructed on lines \textit{BA} and \textit{AC} are all depicted on the accompanying diagram (see Fig. \ref{figI47}). 

Yet, already in the first propositions of Book II, we observe a new phenomenon: figures mentioned in the \emph{diorisomos}  are not represented on the diagrams. Here are \emph{ekthesis} and \emph{diorisomos} of  proposition II.2:
``For let the straight-line AB have been cut, at random,
at point C. I say that the rectangle contained by AB, BC together with the rectangle contained by BA, AC, is equal
to the square on AB." 

On the accompanying diagram (see Fig. \ref{figII2}), we can see the lines \textit{AB}, \textit{AC}, and \textit{BC},  however, neither the rectangle contained by \textit{AB}, \textit{BC}, nor the one contained by\textit{BA}, \textit{AC} is  depicted on the diagram:
 line-segments \textit{AB}, \textit{AC}, and \textit{BC}   lay  on the same straight line and do not contain a right-angle.  
 Rectangles which are supposed to be formed by line segments  not containing a right-angle occur in every proposition of Book II. 

Furthermore, notice  that in proposition II.14, Euclid shows how to ``construct a square equal to a given rectilinear figure".
Although the accompanying diagram depicts a quadrilateral $A$,  it is a symbol of a rectilinear figure  rather than the individual object that is being  studied in the proposition: it is not constructed,  its vertices are not denoted by letters, 
and they are not   involved in the constructions carried out in  the proposition.

 The diagram in proposition VI.31 plays an even more abstract role. The \emph{diorisomos}  part of VI.31 reads:  ``Let ABC be a right-angled triangle having the angle ABC a right-angle. I say that the figure  on BC is equal to similar, and similarly described,
figures on BA and AC." Like  figure \textit{A} in proposition II.14, the rectangles on the sides  of the triangle \textit{ABC} are not constructed,  they are not involved in the proof, and the proof does not rely on information  read-off the diagram.\footnote{For a detailed analysis of this proposition in terms of individual and abstract components, see (Błaszczyk, Petiurenko 2020).}

\subsection{\emph{Rectangle contained by two straight-lines}}

Here is the crucial  definition of Book II: ``Any right-angled parallelogram is said to be contained by two straight-lines containing a right-angle".
Throughout the entire book, Euclid studies squares, rectangles, and triangles. The term \textit{right-angled parallelogram} is applied only to rectangles, then it is reshaped to \textit{rectangle contained by two straight-lines}. What, then, is the role of the term  \textit{rectangle contained by two straight-lines}? How does it differ from a simple rectangle?

Let us start with the most general concept, namely that of a figure. In Book I, Euclid defines a  figure as follows: ``A figure is that which is contained by some boundary  or boundaries". The term \emph{boundary} applies to a circle only, \emph{boundaries} apply to polygons. Hence, for example, a triangle is contained by three straight-lines, i.e., its sides. In other words, what we call  a polygonal curve today is not considered to be a single line in the \emph{Elements} -- according to Euclid,  it is a composition of lines.

In Book II, in addition  to triangles,  Euclid studies squares and rectangles. Definition I.22 clarifies these  concepts. It reads: ``And of quadrilateral figures: a square is that which is right-angled, and equilateral, a rectangle that which is right-angled but not equilateral".  The term \textit{parallelogram} requires the concept of parallels and is not included in the list of definitions prefacing Book I. It occurs in proposition I.34 as  \textit{parallelogrammic figure}. Although  not explicitly defined,  it is clear what Euclid means: a parallelogram  is a quadrilateral with two pairs of parallel sides (I.33 shows the existence of parallelograms).  This term is closely related to Euclid's theory of equal figures. Within this theory, in proposition I.44, Euclid shows how to construct a parallelogram when its two sides and an angle between them are given. Still, in Book II, all parallelograms are rectangles. What is, then, the reason for the term   \emph{rectangle contained by two straight lines}?

Firstly,  this term   is related to the ways figures are referred to in the text of the propositions, specifically,  it is essential in \emph{protasis} parts. 
Secondly, it  plays an  analogous role to the term \textit{square on} a side:  as the  latter  enables to identify a square with one side,  the former enables to identify a rectangle with  two sides with no reference to a diagram. Thirdly, since the terms \textit{square on a straight-line} and \textit{rectangle contained by straight-lines} are applied both to figures which are represented (visible) and not represented (invisible) on a diagram, these names make it possible to relate  the visible and invisible figures. Due to substitution rules which we  detail  in section \S\,5, Euclid can claim that a  \textit{rectangle contained by X,Y}, which is not represented on the diagram, is \textit{contained by A,\,B}, where segments \textit{A},\,\textit{B} form a rectangle which is represented on the diagram. 

%To be clear, these substitution rules apply to the relation \textit{contained by} rather than  an equality of \textit{rectangles contained by}. Usually, a rectangle \textit{contained by} taken together with another figure is equal to a figure.
%In II.11 and II.14, we can find a \textit{rectangle contained by} which is equal to a square represented, and not represented on the diagram respectively. Yet, these equalities result from a combination of the Pythagorean theorem and Common Notion 2. In other words, whenever a \textit{rectangle contained by} is equal to another figure, it is not a straightforward relation.

\subsubsection{Current interpretations of \textit{rectangle contained by}}

Within the so-called geometric algebra interpretation of Book II,  all rectangles  are represented by the  formula $ab$, no matter whether lines $a, b$  contain a right-angle or not. Moreover, the term $ab$ is  subject to some processing.

{To be clear, we do not agree with the claim that this interpretation ignores the historical context by implying a multiplication of line segments. One may treat the  terms  
$ab$, $a(b+c)$ as interpreting the phrases \textit{rectangle contained by} lines  $a, b$, or a rectangle contained by lines $a$ and $b, c$. Yet, when that the term $ab$ is applied in the same way to rectangles represented and not represented on a diagram, it blurs the essence of Book II.}

John Baldwin and Andreas Mueller interpret a \textit{rectangle contained by X, Y} as a rectangle  of length \textit{X} and height \textit{Y}. They write: ``This definition [II.1] allows us to study the areas of arbitrary rectangles by indexing each rectangle by its semi-perimeter and a cut point in that line" (Baldwin, Mueller 2019, 8). 
In proofs, instead of \textit{Y}, they consider \textit{W} such that \textit{W=Y}, and lines \textit{X}, \textit{W} contain a right-angle. Thus, in fact, they reduce a \textit{rectangle contained by} to a rectangle represented on a diagram.  Still, they are the only authors that find Euclid's definition ``strange to modern ears". 
%Both interpretations view  the term \textit{rectangle contained by} from the perspective of a dissection of figures. 

Ian Mueller tries yet another trick: ``\textbf{O}(\textit{x,y}) is used to designate  a rectangle with arbitrary straight lines equal to $x$ and $y$  as adjacent sides" (Mueller 2006, 42). First of all,  $\textbf{O}(x,y)$ denotes a rectangle and makes no difference between visible and invisible figures. For example, in Mueller's formalization of II.2, there is no difference between  \textbf{O}(\textit{AD,AC}) and \textbf{O}(\textit{AB,AC}), since $AB=AC$ (see Fig \ref{figII2}). By his notation alone, Mueller ignores   the basic problem Euclid seeks to resolve in propositions II.1--8, namely the relation between a rectangle represented on the diagram, and  one not represented on the diagram. As a result, he distorts Euclid's original proofs, even though he can easily interpret the theses of his propositions.\footnote{In fact, Mueller tries to reconstruct only the proof of II.4. For  details, see \S\,7 below.}   

Jeffrey Oaks provides a similar interpretation, as he writes in a commentary to proposition VI.16 of the \textit{Elements}:
``Here `the rectangle contained by the means' in most cases will not be a particular rectangle given in position because the two lines determining it are not attached at
one endpoint at a right angle. In fact, the sides determining rectangles cited in Greek
works rarely satisfy Euclid’s definition at the beginning of Book II [...]. Already in Proposition II.1 Euclid writes about `the rectangle contained by
A, BC' when the two lines may not be anywhere near each other. And the lines determining the rectangles cited in Proposition II.2 are absolutely not at right angles,
since they are colinear. Propositions VI.16 and XI.34, like many propositions in Greek
mathematics, are about the measures or sizes of geometric objects. `The rectangle
contained by the means' does not designate a particular rectangle given in position,
but only the size of a rectangle whose sides are equal (we would say “congruent”) to
those lines. Location and orientation of this rectangle relative to the other magnitudes
in the diagram are undetermined and irrelevant to the argument. It is only the relative
`measure' that is intended" (Oaks 2018, 259).

At first, Oaks admits that not all \textit{rectangles contained by}
are featured on diagrams. Going beyond this  observation, he links (invisible) \textit{rectangles} with the concept of measure.

The above interpretations involve terms,  figures  and measure. We do not interpret the phrase \textit{rectangle contained by},  but rather study its role in Euclid's proofs. Eventually we view it as a proof technique not an object.

\subsection{Identifying figures through letters and technical terms}

In the \emph{Elements}, Euclid adopts the following pattern of naming the figures featured on diagrams: squares and rectangles  are, first of all, denoted by the letters located next to their vertices, they are also denoted by the letters which  designate the diagonal.
In proposition I.46, Euclid shows how to describe a square on a given straight-line. In the propositions that follow, squares are also identified by the  phrase \textit{square on a straight-line}, where the  specific name of a line is given.
We  can illustrate this naming technique  by referring to proposition I.47 (Fig. \ref{figI47} represents the accompanying diagram).

Thus, in the text of the proposition, the square \textit{BDEC} is also called \emph{the square on BC}; the square on \textit{BA} is also denoted by the two letters located on the diagonal, namely  \textit{GB}.  Since the intersection of lines \textit{BC} and \textit{AL} is not named, rectangles that make up  the square \textit{BDEC} are named with two letters, as \emph{parallelogram BL} and \emph{parallelogram CL}.

 In Book II,  Euclid introduces yet another  naming scheme for the rectangle: it is identified with its two sides and is called  \textit{rectangle contained by}, and the term is followed by the names of line-segments containing a right-angle. 
 
\begin{figure}
\begin{center}
\includegraphics[scale=1.2]{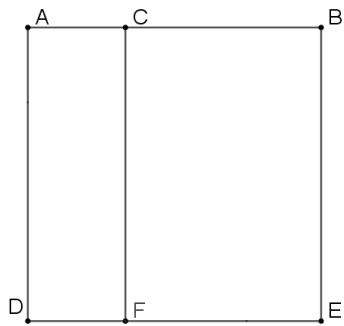}
\caption{\textit{Elements}, II.2.} \label{figII2}
\end{center}
\end{figure}

 All these naming rules -- that is,  by  vertices,  by a diagonal, \textit{square on}, and \textit{rectangle contained by} --  apply to figures represented on the diagrams.  Here, for example, is the text of proposition II.2  (Fig. \ref{figII2} represents  the  accompanying diagram). 
 
 [1] ``For let the straight-line AB have been cut, at random,
at point C. I say that the rectangle contained by AB, BC together with the rectangle contained by BA, AC, is equal
to the square on AB. 

[2] For let the square ADEB have been described on AB,  and let CF have been drawn through C
parallel to either of AD or BE.

[3] So  AE is equal to the  AF, CE. And AE is the square on AB. And AF (is) the rectangle contained by 
BA, AC.  For it is contained by DA,  AC, and AD (is) equal to AB. And CE (is) by AB, BC. For BE (is) equal to AB.
AB. Thus,  by BA, AC, together with  by AB, BC, is equal to the square on AB."

Thus, the square \textit{ADEB} is also named \textit{the square on AB}.
Rectangle {AF} is also called \textit{the rectangle contained by DA,  AC}.

The phrases \textit{square on} and \textit{rectangle contained by} are also applied to figures not represented on  diagram. In the text of proposition II.4,  the term \textit{the square on AC} occurs, although there is no such  square on the accompanying diagram (see Fig. \ref{figII3_4}). Also in this proposition, the term 
\textit{rectangle contained by AC, CB} occurs, although there is no such  rectangle on the accompanying diagram.

\begin{figure}
\begin{center}
\includegraphics{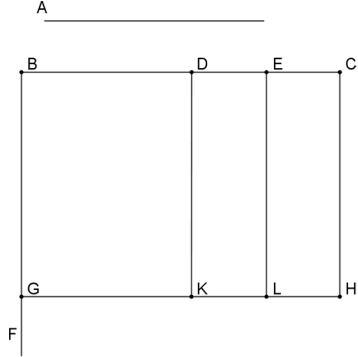}
\caption{\textit{Elements}, II.1.} \label{figII1}
\end{center}
\end{figure}

In fact, rectangles contained by straight-lines lying on the same line and not containing a right-angle are common in Book II. Whether represented  on  diagrams or not, as a rule, they are contained by individual straight-lines. However,  proposition II.1 represents a unique case in this respect.  Therein, Euclid considers rectangles contained by \textit{A},\,\textit{BD}, and 
\textit{A},\,\textit{DE}, and \textit{A},\,\textit{EC} (see Fig. \ref{figII1}). They are to 
be   rectangles contained by   \textit{BG}, \textit{BD},   by  \textit{DK}, \textit{DE}, and  by \textit{EL}, \textit{EC} respectively. 
 Rectangles contained by \textit{A}, \textit{BD}, by \textit{A}, \textit{DE}, and \textit{A}, \textit{EC} are neither represented on the diagram, nor contained by individual line-segments: line \textit{A},  considered as a side of these rectangles,  is not an individual line.
 
 To be clear, the relation between rectangles,  contained on the one hand by \textit{A},\,\textit{BD}, and on the other  by \textit{BG},\,\textit{BD}, is not an equality: it is stated that the former is the rectangle contained by \textit{BG},\,\textit{BD}. The relevant part of the proposition reads: ``BK is by A, BD. For it is contained by GB,\,BD, and BG is equal to A".  Here, \textit{BK} is represented on the diagram, and Euclid  claims that it is  contained by \textit{BG},\,\textit{BD}, which is simply another name of the rectangle \textit{BK}. Still, Euclid also claims that \textit{BK} is contained by \textit{A},\,\textit{BD}, while the later rectangle  is not represented on the diagram.  Thus, the relation between figures represented and not represented on a diagram is founded on  substitution
 and renaming  rules.  We will explicate these rules in  section \S\,5. In general, Euclid engages a bunch of  
  tricks to establish an equation between a \textit{rectangles contained by} and a figures represented on a diagram.

Interestingly, Euclid never refers to  proposition II.1. Moreover, its result, when viewed from the modern perspective, is reiterated in propositions II.2 and II.3. Hence, it seems that  its role  is to demonstrate the substitution rules  which are applied throughout the rest of Book II, rather than to present a specific geometrical statement. All that is required to analyze these rules are propositions II.1--II.4. 
From the perspective of substitution rules, proposition II.1 introduces them, then proposition II.2 applies them to  \textit{rectangles contained by}, and proposition II.4 -- to \textit{squares on}. Proposition II.4 involves yet another object, namely the so-called \textit{complement}. It shows how to apply the  substitution rules to these objects. 

From the perspective of represented vs not represented figures, proposition II.2 equates figures which are represented, on  the one side, and not represented, on the other, while
proposition II.3 equates figure not represented, on the one side, and figures represented and not represented, on the other side, proposition II.4 introduces yet another operation on figures which are not represented, as it includes an object called \textit{twice rectangle contained by}, where the rectangle is not represented on the diagram.

\subsubsection{Algebraic view on propositions II.1--3}
 Without paying attention to Euclid's vocabulary, specifically to the terms  \textit{square on} and \textit{rectangle contained by}, one cannot find a reason for propositions II.2 and II.3. Thus, Bartel van der Waerden in (Waerden 1961) considers them
 as special cases of II.1. Similarly, 
Robin Hartshorne, in the \textit{Appendix} to (Hartshorne 2000), includes statements of ``the most frequently quoted results" of the \textit{Elements}. Regarding Book II, he refers to proposition II.1, then skips to II.4.  
%Richard Fitzpatrick in his translation of the \textit{Elements} (Fitzpatrick 2007), adopts alike interpretation.

  From the modern perspective, especially when the \textit{diorismos} of Euclid's proposition is 
stylized as an algebraic formula, such an interpretation seems reasonable. For, when II.1 is rendered as $a(b+c+d)=ab+ac+ad$, then II.2 is $a^2=ab+ac$, given $a=b+c$, and II.3 is $a(b+a)=ab+a^2$. Indeed, II.2 and II.3 follow from II.1 by suitable substitutions.  In fact, however,  proposition II.1 is never quoted in the \textit{Elements}, and 
due to the role of line \textit{A} it is a unique proposition in the entire \textit{Elements}.

  \section{Schemes of propositions II.1 to II.4}

In this section, we provide detailed analysis of propositions II.1 to II.4. They aim to reveal non-geometrical rules
which enable to relate figures represented and not represented on the diagrams.

  Here is  the text of proposition II.1, starting with the \textit{diorismos}\footnote{Numbering of sentences and  names of parts  added.}:

\textit{Diorismos} ``Let A, BC be the two straight-lines, and let BC, be cut, at random,
at points D,\,E. I say that the rectangle contained by A, BC is equal to the rectangles contained by
A, BD, by A, DE, and, finally, by A, EC."

\textit{Kataskeu\={e}}\ \ ``For let BF have been drawn from point B, at right angles to BC, and let BG be made 
equal to A [...]."

 \pvn``[1] So   BH is equal to  BK, DL, EH.
 [2] And BH is by A, BC. For it is contained by GB, BC, and BG is equal to A.
  [3] And BK (is)  by A, BD. For it is contained by GB, BD, and BG is equal to A.
     [4] And DL is by A, DE.  For DK, that is to say BG, is equal to A.
[5] Similarly, EH (is) also by A, EC.
[6] Thus, by A, BC is equal to by A, BD, by A, DE, and, finally, by A, EC."

Now, we present this proposition in a more schematic form. In what follows, symbol  $A\times BC$ stands for the phrase ``rectangle contained by A,\,BC", while
 $BH\, \pi\, GB\times BC$ stands for the phrase ``BH is contained by GB,\,BC".
 
\pvn\textit{Diorismos}
\[A\times BC=A\times BD,\, A\times DE,\, A\times EC.\]
\textit{Kataskeu\={e}} $BF\perp BC,\ BG=A$; $GH\parallel BC$; $ DK, EL, EH\parallel BG$.\\
\textit{Apodeixis}
\begin{eqnarray*}
{\color{red}BH=BK, DL, EH}  && \\
 {\color{blue} BH\, \pi\, GB\times BC},\  BG=A              &\xrightarrow[]{}& {\color{violet}BH\,\pi\,A\times BC} \\
{\color{blue}BK\,\pi\, GB\times BD},\  BG= A&\xrightarrow[]{}& {\color{violet}BK\,\pi\, A\times  BD}\\
    DK=BG=A               &\xrightarrow[]{}& {\color{violet}DL\,\pi\,A\times DE}\\
    &\xrightarrow[]{}& {\color{violet}EH\,\pi\,A\times EC}\\
  &\xrightarrow[]{}& {\color{magenta}A\times BC= A\times BD,\, A\times DE,\, A\times EC}.
\end{eqnarray*}

(1) The  formula in red interprets sentence [1].  It is a simple statement with no justification and  the starting point of the whole argument.  Since the figures involved are represented on the diagram, we  interpret it as based on  purely visual evidence.\footnote{(Błaszczyk, Petiurenko 2020) develops the idea of pure visual evidence.} 

All of the rectangles mentioned in the \textit{diorismos} are not represented on the diagram. 
Therefore, we are to explain how, starting from the relation between the figures represented on the diagram, Euclid gets the relation between figures not represented on the diagram.

(2) The next line in the \textit{apodeixis} scheme interprets sentence [2]. The formula in blue stands for the phrase
  ``BH [...]  is contained by GB, BC". It is one of the three possible names for a rectangle represented on a diagram. 
  Thus, it is a result of renaming figures rather than a geometrical or logical relation. Hereafter, formulas resulting from renaming will be represented in blue.
  
 The equality $BG=A$ follows from the construction part of the proposition.  Yet, the most puzzling is the phrase ``BH is by A, BC". We interpret it as a result of substituting \textit{A} to the formula $BH\,\pi\, BG\times BC$ in place of \textit{BG}. Arguments of this kind are  applied all throughout Book II.     The relation between the rectangle \textit{BH}, and the one contained by \textit{A},\,\textit{BC}  is by no means an equality; the word pattern ``BH is by A, BC" is systematically used by Euclid.
   Let us represent  formulas obtained by this type of substitution in violet.

(3)  In sentence [3], Euclid reiterates the previous argument. %Therefore line (3) in our scheme does not introduce any new idea.  

(4) In sentence [4], Euclid skips supposition ${\color{blue}DL\,\pi\,DK\times DE}$ and notes the equalities 
 $DK=BG=A$. They can be justified by \mbox{Common Notion 1.}
 
 (5) In sentence [5], Euclid skips arguments relying on substitution and CN1, and  simply states  the result. It is a  way of shortening repeated arguments, typical of Euclid. 
 
 (6) In sentence [6],  with \foreignlanguage{polutonikogreek}{>'ara},  Euclid reaches the equality between invisible rectangles. We interpret this step as a result of another substitution rule: in the  equality  starting in \textit{apodeixis}, $A\times BC$ is substituted for \textit{BH}, then $A\times BD$ is substituted for \textit{BK}, etc.
 
  Let us represent the equality obtained by this type of substitution, i.e., substitutions for equality, in magenta.

The arrow $\rightarrow$ in the scheme stands for a conjunction, usually it is \foreignlanguage{polutonikogreek}{g'ar}. It is by no means suggested to be  a logical implication.

II.2 (see Fig. \ref{figII2})

 In the below scheme, the term $AB^2$ stands for the phrase ``the square on AB". Thus, for example, $AE\ is\ AB^2$ interprets the phrase ``AE is the square on AB", or ``AE is on AB".

\pvn\textit{Diorismos}
\[AB\times BC, BA\times AC =AB^2.\]
\textit{Apodeixis}
\begin{eqnarray*}
{\color{red}AE=AF,CE}  && \\
     {\color{blue} AE\ is\ AB^2}       && \\
    {\color{blue}AF\,\pi\,DA\times AC},\ AD=AB &\xrightarrow[]{}& {\color{violet}AF\,\pi\,BA\times AC}  \\
    BE=AB              &\xrightarrow[]{}& {\color{violet}CE\,\pi\,AB\times BC}\\
    &\xrightarrow[]{}& {\color{magenta}BA\times AC, AB\times BC =AB^2}.
  \end{eqnarray*}

Lines (3) and (4) of the \textit{apodeixis} scheme, represent, again,  typical of Euclid way of shortening repeated arguments: in line (4), Euclid skips the premise \textit{CE is between CF, BC}. 

In this proposition, the \textit{diorismos} equates the figure represented on the diagram, 
that is, the square on \textit{AB}, with figures not represented on the diagram, namely rectangles contained by \textit{AB}, \textit{BC}, and \textit{AB}, \textit{AC}.

\begin{figure}
\begin{center}
\includegraphics{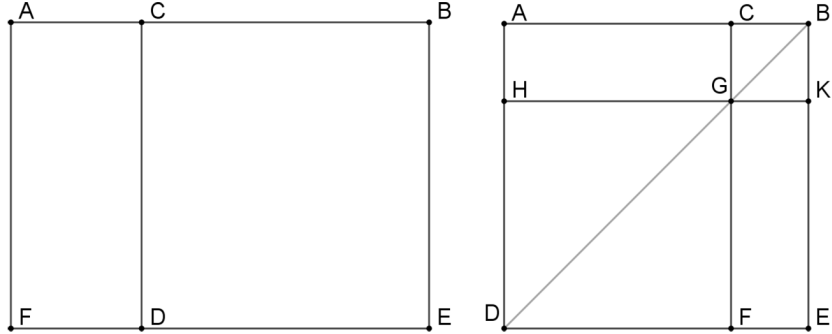}
\caption{\textit{Elements} II.3 (left) and II.4 (right).} \label{figII3_4}
\end{center}
\end{figure}

II.3 (see Fig. \ref{figII3_4})
\pvn\textit{Diorismos}
\[AB\times BC= AC\times CB,\, BC^2.\]
\textit{Apodeixis}
\begin{eqnarray*}
{\color{red}AE=AD,\, CE} && \\
 {\color{blue}AE\,\pi\,AB\times BE},\  BE= BC  &\xrightarrow[]{}&  {\color{violet}AE\,\pi\,AB\times BC}\\
DC=CB    &\xrightarrow[]{}& {\color{violet}AD\,\pi\,AC\times CB} \\
           {\color{blue}DB\ is\ CB^2}
           &\xrightarrow[]{}& {\color{magenta}AB\times BC= AC\times CB,\, BC^2}.
    \end{eqnarray*}

Here, the square \textit{CE} is also named by the second diagonal, as \textit{DB}. Apart from this fact, this \textit{apodeixis} is 
similar to the previous one.

In this proposition, Euclid equates the figure not represented on the diagram, $AB\times AC$, with figures both represented and not represented on the diagram, $AC\times CB,\, BC^2$.

Proposition II.1--3 share the same word patterns clearly represented by the colors of terms in our schemes. 
Since there are no explicit references to these propositions in the rest of Book II, we treat them as introducing a technique of dealing with figures which are not represented on  diagrams.
 
With II.4, we pass to a proposition which will be referred to in the following propositions. 

II.4 (see Fig. \ref{figII3_4})
\pvn\textit{Diorismos}
\[AB^2=AC^2,\, CB^2,\, 2 AC\times CB.\]
\textit{Apodeixis}.
\begin{eqnarray*}
{\color{blue}CGKB\ is\ CB^2} && \\
{\color{blue}HF\ is\ HG^2}, {\color{violet}HF\ is\ AC^2} &\xrightarrow[]{}& HF,\, KC\ are\ AC^2,\, CB^2\\
GC=CB &\xrightarrow[]{}& {\color{violet}AG\,\pi\,AC\times CB}\\
 AG=GE &\xrightarrow[]{}& {\color{magenta}GE=AC\times CB} \\
 &\xrightarrow[]{}& AG, GE=2 AC\times CB\\
HF, CK\ are\ AC^2, CB^2 &\xrightarrow[]{}&  HF,\, CK,\, AG,\, GE=\\
&&=AC^2,\, BC^2,\, 2AC\times CB \\
 {\color{red}HF,\, CK,\, AG,\, CE=ADEB}, &&\\
 {\color{blue}ADEB\ is\ AB^2} &\xrightarrow[]{}& {\color{magenta}AB^2=AC^2,\, CB^2,\, 2 AC\times CB}.
    \end{eqnarray*}
    
    In this proposition, Euclid equates the figure  represented on the diagram, $AB^2$, with figures represented, $CB^2$, and not represented on the diagram, $AC^2$, $2AC\times CB$. A new component is the square on \textit{AC} which is not represented on the diagram.
    Next, ``twice the rectangle contained by AC, CB" means that Euclid  handles not represented rectangles the same way as  represented ones. 

Formula $HF\ is\ HG^2, HG^2\ is\ AC^2$ interprets the following phrase: ``HF is also a square. And it is on HG, that is to say AC". We formalize it like this:  
$$HF\ is\ HG^2, HF=AC \rightarrow HF\ is\ AC^2.$$ 

Thus, from the adopted perspective, it is a substitution rule applied to \textit{the square on}, analogous to the one applied to \textit{the rectangle  contained by}. That is why it is highlighted in violet.

Now, let us focus on the lines (3) and (4) of the \textit{apodeixis} scheme. They interpret the following part of Euclid's proof: ``AG is equal to GE. And AG is contained by AC, CB. For GC is equal to CB. Thus, GE is also equal to the one by AC, CB". The equality $AG=GE$ follows from proposition I.43. $AG\,\pi\,AC\times CB$ is the result of the substitution rule we  already identified. The conclusion ``GE is also equal to the one by AC, CB" is 
the result of a substitution to the equality.

In regard to the term $AG, GE=2 AC\times CB$, we cannot provide a clear justification for this conclusion. The same applies to the conclusion 
\[HF,\, CK,\, AG,\, GE= AC^2,\, BC^2,\, 2AC\times CB.\]

Although we could justify it by using  logical tricks, it is not  Euclid's style to conceal complicated rules and explicate simple ones.  It seems to be a puzzle Euclid could not resolve. 

The formula in red  interprets the following sentence: ``HF, CK, AG, GE are the whole ADEB". Like in previous cases, we take it to be justified by  visual evidence.  
From a logical point of view, it could be the starting point of this proof.

\subsubsection{David Fowler on proposition II.4}

In regard to propositions II.1--8, Fowler writes: ``In all of these propositions, almost all of the text  of the demonstrations concerns the construction of the figures, while the substantive content of each enunciation is merely  read off from constructed figures, at the end of the proof, as in lines 49 to 50 of II.4, the proposition just considered" (Fowler 2003, 69)

Our schemes evidence that arguments read off the diagrams are staring points of the proofs
II.1--3. In II.4--8, they are at the end of proofs. However, the enunciations of propositions II.1--8 concern figures which are not represented on the diagrams, therefore the essential arguments can not be read off the diagrams.

 Commenting on the final lines of the proof II.4, Fowler writes: 

\pvn \textit{Therefore the four areas HF, CK, AG, GE are equal to the square on AC, CB, and twice the rectangle by AC, CB},
\pvn and the fact that this gives a decomposition of the square, the ostensible point of the proposition, is merely stated (lines 49 to 50): 
\pvn \textit{But HF, CK, AG, GE are the whole ADEB, which is the square on AB}" (Fowler 2003, 69)

The point is that \textit{HF}, \textit{CK}, \textit{AG}, \textit{GE} give the decomposition of the square \textit{ADEB}, not  ``AC, CB, and twice the rectangle by AC, CB".  

\subsubsection{Ian Mueller's interpretation again}
Interestingly,   the equality $GE=AC\times CB$ follows from the equality 
$AG=GE$ rather than from the rule: 
$$A\,\pi\, C\times D,\ C=E,\ D=F,\ B\,\pi\, E\times F \rightarrow  A=B.$$ 
 
 This 
alternative argument is obvious within the so-called geometric algebra, as well as Mueller's interpretation (see \S\,3.2 above). More specifically, within Mueller's interpretation the equality (congruence) obtains
\[x=A, y=B\Rightarrow \textbf{O}(x,y)\simeq \textbf{O}(A,B).\] 

However, we have not identified this rule in Book II. On the contrary, proposition II.4 exemplifies different reason, namely
\[A\,\pi\,C\times D,\  A=B \rightarrow B=C\times D.\]

In words, since we know that figures \textit{A} and \textit{B} are equal and one of them is \textit{contained by} \textit{C},\,\textit{D}, then the other is also \textit{contained by} \textit{C},\,\textit{D}. It means, that the equality established within the theory of equal figures is more fundamental than the relation \textit{contained by}. In a way, Euclid aims to  introduce \textit{rectangles contained by} into equalities of non-congruent figures. 
  Within Mueller's interpretation the equality between rectangles follows from a relation \textit{contained by}.

\section{Non-geometrical rules in Book II}

Four colors in or schemes of Euclid's propositions correspond to three groups of rules: visual evidence (red), renaming (blue) and substitutions (violet and magenta). In this section we scrutinize these rules.
\subsection{Visual evidence}

It is  standard to identify    two meanings of equality of figures in the \textit{Elements} --  congruence and the equality of non-congruent figures. 
The congruence of figures is usually linked to the idea of coinciding figures  involved in \textit{Common Notion} 4. It is also commonly assumed that the idea of coinciding figures  plays a crucial role in proposition I.4. But the superposition  of figures presupposes  (rigid) motion. 

The statements in red in our schemes are so simple that they do not engage any other concepts. If any needed justification,  CN4 would be a good choice. Although such an interpretation finds no textual corroboration, there are no significant differences between the claim that CN4 is founded on visual evidence and the claim that statements in red are justified by CN4. In fact, Euclid provides no arguments for these statements.

\subsection{Overlapping figures}

At first, the range of visual evidence is obvious: it justifies the equality of a figure and their dissection parts, which are squares and rectangles.  Yet, as we proceed further, in propositions II.5--8,  Euclid extends its power to another cases. In II.5, the gnomon \textit{NOP} is taken together with the square \textit{LG} and Euclid declares that they form the square \textit{CEFG}: ``the gnomon NOP and the square LG are the whole square CEFB" (see Fig. \ref{figII5}). We interpret it as an equality based on visual evidence, $NOP, LG=CEFB$.  

In the next  proposition, the diagram is the same as regards the gnomon and its complementing square. This time, instead of the square \textit{LG}, Euclid adds  the square on \textit{BC}  to the gnomon \textit{NOP} (see Fig. \ref{figII6}). His argument is this:
since $LG=BC^2$, then $NOP, BC^2=CEFB$. However, while the equality $NOP, LG=CEFD$ is visually obvious, the equality  $NOP, BC^2=CEFB$ requires other kind of justification, as $BC^2$ is not represented on the diagram.
 If we add an auxiliary line to represent the square on \textit{BC}, we will get overlapping figures. Thus, the status of the equality $NOP, BC^2=CEFB$ is not as obvious as $NOP, LG=CEFB$ in II.5.
 
\begin{figure}
\begin{center}
\includegraphics{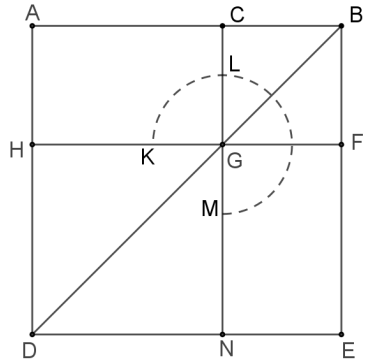}
\caption{\textit{Elements} II.7.} \label{figII7}
\end{center}
\end{figure}

 In II.7, Euclid goes a bit further in terms of abstraction. He explicitly considers  overlapping figures: the rectangle  \textit{AF}, and the square \textit{CE}. He claims that they together form the gnomon \textit{KLM}  and the square \textit{CF}: ``but AF, CE are the gnomon KLM and the square CF" (see Fig. \ref{figII7}). Here is the scheme of the proof.
%\subsection{Overlapping figures}

\pvn\textit{Diorismos}
\[AB^2, BC^2= 2 AB\times BC, CA^2.\]
\textit{Apodeixis}

\begin{eqnarray*}
AG=GE &\xrightarrow[CN3]{}& AG, CF=GE, CF \\
 AF=CE&\xrightarrow[]{}&  AF, CE= 2 AF\\
  &\xrightarrow[]{}& {\color{red}KLM, CF=2 AF}\\
  BF=CB &\xrightarrow[]{}&  {\color{magenta}2AF= 2 AB\times BC}\\
 {\color{blue} DG\ is\ AC^2} &\xrightarrow[]{}& KLM, DG, GB= 2AB\times BG, AC^2\\ 
 {\color{red}KLM, DG, GB=ADEB, GB} &\xrightarrow[]{}& ADEB, GB = 2AB\times BG, AC^2\\
  &\xrightarrow[]{}& AB^2, BC^2=  2AB\times BG, AC^2.   
    \end{eqnarray*}
    
    Here,  the red suggest that the respective formulas are  based on  visual evidence. Yet, since rectangles \textit{AF} and \textit{CE} overlap, i.e., share the square \textit{CF}, they do not represent the same kind of evidence as red formulas  in propositions II.1--4.
    
\begin{figure}
\begin{center}
\includegraphics{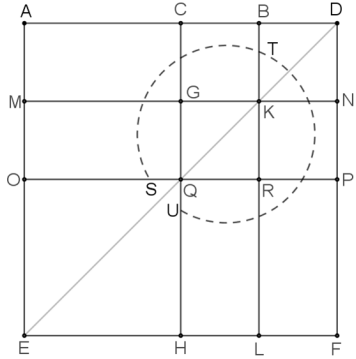}
\caption{\textit{Elements} II.8.} \label{figII8}
\end{center}
\end{figure}

    In II.8, Euclid's considers even more complicated configuration of overlapping figures. Below we present the scheme of the key step (see Fig. \ref{figII8}).
    
    \begin{eqnarray*}
DK=CK=GR=RN &\xrightarrow[]{}& DK, CK, GR, RN=4 CK \\
 AG=MQ=QL=RF&\xrightarrow[]{}&  AG, MQ, QL, RF= 4 AG\\
    &\xrightarrow[]{}& {\color{red}DK, CK, GR, RN,AG, MQ, QL, RF=STU}\\ 
    &\xrightarrow[]{}& {STU=4 AK}.   
    \end{eqnarray*}
 
 The red formula interprets the following sentence: ``the eight, which comprise the gnomon STU". 
Due to an implicit step  
$$DK, CK, GR, RN, AG, MQ, QL, RF=4 AK,$$ 
Euclid managed to bypass an explicit reference to overlapping figures. The missing step  could be like this: \textit{AG=MQ=QL}, \textit{CK=GR}, then \textit{AK=MR=GL}.  However, here the square \textit{GR} is counted twice. Thus,  the square \textit{GR} could be moved to cover the square \textit{DK}. Argument of this kind characterize the so-called  dissection proofs, for example, the famous Chinese proof the Pythagorean theorem.     
 Significantly, Euclid does not apply dissection combined with a translation. 
 
 Another option could be a reference to CN2, namely: since \textit{AG=QL} and \textit{CK=DK}, then \textit{AG, CK= QL, DK}. However, as a rule, Euclid does not apply CN2 when the resulting figure  is not connected, i.e., it does not make a \textit{whole}. Nevertheless, it is possible when the resulting figures overlap. 
 
 With the use of modern technology, overlapping figures are handled with shades, or colors. 
Yet,  these are textbooks tricks.  Foundational studies seek to eliminate overlapping figures, as we demonstrate in the next section.

\subsubsection{Hilbert-style account of visual evidence}
Within Hilbert's tradition of reading  the \textit{Elements}, the congruence of line segments, angles, and triangles is 
covered by their respective axioms,  Euclid's proposition I.4 specifically  is an axiom in the Hilbert system.  Euclid's theory of equal figures is covered by the idea of the content of a figure. 

Robin Hartshorne develops Hilbert's idea of content further to the modern concept of measure. Regarding Euclid's theory of equal figures, he writes: ``Looking at Euclid's theory of area in Books I--IV, Hilbert saw how to give it a solid foundation. We define
a notion of equal content by saying that two figures have equal content if we transform one figure into the other by adding and subtracting congruent  triangles" (Hartshorne 2000, 195).

Indeed, figures involved in   equalities in red  also have equal content in the  Hilbert-Hartshorne system.
 However, a justification is far from obvious. Let us take, for example, proposition II.2 and Euclid's statements \textit{AE=AF,\,CE}. Hartshorne's definition of a figure reads: ``A rectilinear figure [...] is a subset of the plane that can be expressed as a finite non-overlapping union of triangles" (Hartshorne 2000,  196).   \textit{AE}, on the one hand, and  \textit{AF}, \textit{CE}, on the other,  meet  the requirements of this definition. The next  definition is this:
``Two figures $P, P'$ are equidecomposable if it possible to write them as non overlapping unions of triangles
\[P=T_1\cup....\cup T_n,\ \ P'=T'_1\cup....\cup T'_n, \]
where for each $i$, the triangle $T_i$ is congruent to the triangle $T'_i$" (Hartshorne 2000,  197).
Then, the definition of equal content follows: ``Two figures $P, P'$  have equal content if there are other figure $Q, Q'$ such that: (1) $P$ and $Q$ are not overlapping, (2) $P'$ and $Q'$ are not overlapping, (3) $Q$ and $Q'$ are equidecomposable, (4) $P\cup Q$ and $P'\cup Q'$ are equidecomposable" (Hartshorne 2000,  197).

Thus, the proof that \textit{AE} and \textit{AF,\,CE} have equal content would be the same as the proof of the reflexibility of the relation \textit{have equal content}, as if \textit{AE} and \textit{AF,\,CE} were the same figures. In fact, from the perspective of  set theory, which makes  the basis of the Hartshorne system, $AE=AF\cup CE$.  However, it is not enough. To decide that \textit{AE} and \textit{AF,\,CE} are equidecomposable, we not only have to cover both sides with the same  triangles, we also need to add to both sides another figure \textit{Q}. This peculiar step  is the price for the \textit{solid}  account of Euclid's visual evidence.

Now, let us take the rectangle contained be \textit{BA,\,AC}. Given Hartshorne's definition, it is not a figure at all. Therefore $BA\times AC$  cannot be studied within this theory.

In sum, within the Hartshorne system,  one can provide conceptually  complicated proof of a  statement which is obvious in the \textit{Elements}. However, a complete reconstruction  of Book II  is impossible since there is no counterpart of the concept \textit{rectangle contained by}.  Hartshorne overestimates his system when he claims that   ``In Book II, all of the results make statements about certain figures having equal content to certain others, and all of these are valid in our framework" (Hartshorne 2000, 203)  In fact, his system   does not enable  to identify the real problems of Book II, that is, a relation between the represented and not represented figures. 

In modern system of geometry,  the measure of a figure, that is a real number, plays  the role of figures which are not represented.

\subsection{Renaming}

 Our schemes of Euclid's propositions clearly expose the role of the names of figures in the analyzed arguments.
Rectangles represented on the diagrams are named by their vertices,  diagonals, and as \textit{contained by} two line segments.
Squares represented on the diagrams, similarly,  are named by their vertices, diagonals, and as  a \textit{square on} a side. Figures which are not represented get only one name: it could be a \textit{rectangle contained} by two lines, 
or a \textit{square} on a line.
Thus, the most important factor is that figures represented on the diagrams can also be named \textit{rectangle contained by}, or \textit{square on}.  Then, due to substitution rules, they can be related to figures not represented on the  diagrams.   

In a model example, in proposition II.2 (see Fig. \ref{figII2}), the rectangle \textit{AF} is represented on the diagram and gets the name  \textit{contained by DA,\,AC}. Segments \textit{DA,\,AC} are represented on the diagram and contain the right-angle. 
Then, Euclid claims that ``AF is contained by BA, AC", for ``AD is equal to AB". However, the rectangle contained by \textit{BA, AC} is not represented on the diagram. Moreover,  these lines do not contain a right-angle. That is why, in our scheme,  $AF\,\pi\,DA\times AC$ is represented in blue -- it is simply a new name for a visible figure. Nevertheless, to turn $AF\,\pi\, DA\times AC$ into $AF\,\pi\,BA\times AC$  a substitution rule is needed, namely rule (3) presented in the next section.

\subsection{Substitution}
First of all, observe that it is not explicit that the relation  \textit{contained by} is  commutative, therefore  we will not apply the following rule
$X\times Y=Y\times X$.\footnote{Mueller writes: ``Since Euclid normally takes for granted such geometrically
obvious assertions as $\textbf{T}(x)\simeq \textbf{O}(x,x)$ and $\textbf{O}(x,y)\simeq \textbf{O}(y,x)$, he could have carried out geometrical versions of theses arguments" (Mueller 2006, 46). However, we have not identified such arguments in Book II.}   Euclid also does not apply this seemingly obvious rule: if $X=U, Y=W$, then $X\times Y= U\times W$.

The first substitution rule, applied all throughout  Book II, is the one denoted in violet in our schemes. It is as follows
\begin{equation}X\,\pi\,Y\times V,\ Y=U\Rightarrow X\,\pi\,U\times V.\end{equation}

The point is, while $X$ and $Y\times V$, as well as $Y$, $V$, and $U$, are represented on the diagram, $U\times V$ is not.
The following line from the scheme of proposition II.2 exemplifies this rule:\footnote{It often happens that Euclid permutes letters naming line segments, as here with \textit{AB} and \textit{BA}, or \textit{AD} and \textit{DA}. This could be a topic for another paper. Generally, it seems that these letters are arranged to follow the drawing of the line, which is to illustrate an argument.}
\[AF\,\pi\, DA\times AC,\ AD=AB \rightarrow AF\,\pi\,AB\times AC.\]

A similar rule applies to \textit{square on}, namely
\begin{equation}  X=Y^2, Y=U\Rightarrow X=U^2.    \end{equation}

The point is, while the square $X$ and its side $Y$ are represented on a diagram, the square $U^2$ is not, although the side $U$ is represented.

We exemplify it by an argument from proposition II.4:
\[HF\ is\ HG^2,\ HG=AC\rightarrow HF\ is\ AC^2.\]

This formula interprets the following phrase: ``HF is also a square. And it is on HG, that is to say on AC".
 Results based on this rule could be also achieved by a reference to proposition I.36. Yet, it would require introducing another point and an extra construction. Significantly,  
Euclid does not refer to I.36 in this context.

Finally, the rule concerning substitution to an equality; in our schemes it is represented in magenta. It is as follows:

\begin{equation}X=Y, \ X\,\pi\, U\times W \Rightarrow U\times W = Y.
\end{equation}

Since the relation of equality is symmetric, by applying this rule, we can also get the following result
\[X=Y,\ \ X\,\pi\, U\times W,\ Y\,\pi\, Z\times V \Rightarrow U\times W = Z\times V. \]

Thus, in proposition II.1, the starting point is this
\[BH=BK, DL, EH.\]

Then, by  rule (1), we get the following results
\[BH\,\pi\,A\times BC,\ BK\,\pi\,A\times BD,\ DL\,\pi\,A\times DE,\ EH\,\pi\,A\times EC.\]

Finally, by  rule (3), we reach the conclusion
\[A\times BC=A\times BD,\ A\times DE,\ A\times EC. \]

To be clear, these substitution rules apply to the relation \textit{contained by} rather than  an equality of \textit{rectangles contained by}. 
In II.11 and II.14, we can find a \textit{rectangle contained by}   equal to a square represented on the diagram.
The equality is achieved by rule (3). Another way to equate \textit{rectangle contained by} and a figure represented on a diagram is by combining the Pythagorean theorem and Common Notion 2. This trick is applied in
 II.11--14.  In other words, whenever a \textit{rectangle contained by} is equal to another figure, it is not a straightforward relation.

\section{From visible to invisible and backward}

In this section, we study the use of propositions II.5,\,6 in II.11,\,14, and  show how through the technique of \textit{rectangles contained by} Euclid  has managed to establish a relation between visible figures.  From a methodological point of view, he applies results obtained in one domain to determine results in another domain. 
It is like a factorization of real polynomial by its factorization in the domain of complex numbers, or,  finding  a solution to a problem in the domain of hyperreals, then, with its standard part, going back to  the domain of real numbers.

\subsection{Propositions II.5--6}

Propositions II.5, 6 are often discussed in the literature, as  scholars seek to provide a reason for including these seemingly twin propositions  into Book II.
First, we show how the  substitution rules   impact the interpretation of these propositions.

\begin{figure}
\begin{center}
\includegraphics{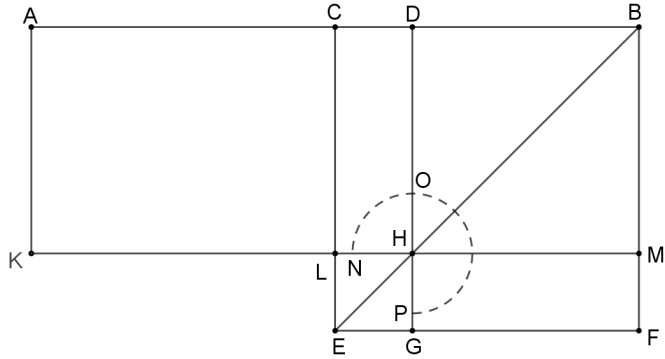}
\caption{\textit{Elements} II.5.} \label{figII5}
\end{center}
\end{figure}

\begin{figure}
\begin{center}
\includegraphics{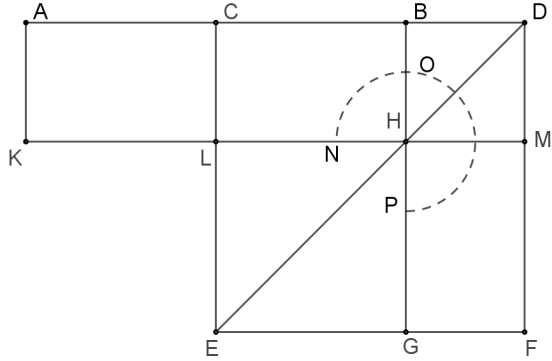}
\caption{\textit{Elements} II.6.} \label{figII6}
\end{center}
\end{figure}

Viewed in terms of construction, they look alike (see Fig. \ref{figII5} and \ref{figII6}). Line \textit{AB} is cut in half at \textit{C}, then point \textit{D} is placed between \textit{C} and \textit{B}, or on the prolongation of \textit{AB}. Yet, their \textit{protasis} parts  differ in wording: in the first case, Euclid considers  \textit{equal} and \textit{unequal lines}, in the second case, \textit{the whole line} and  \textit{the added line}.   Still, when we proceed to their  diagrams and \textit{diorismos},   they are again similar. Moreover,
their proofs apply the same trick: at first, Euclid shows that a rectangle is equal to a gnomom, then he adds  a square that complements the gnomon to a bigger square.

We present a scheme  of proposition II.5  starting from  when it is established  that  the rectangle \textit{AH} is equal to the gnomon \textit{NOP}.

II.5
\pvn\textit{Diorismos}
\[AD\times DB, CD^2=CB^2.\]
\textit{Apodeixis}
\begin{eqnarray*}
... &\xrightarrow[]{}& AH=NOP\\
DH=DB &\xrightarrow[]{}& {\color{violet}AH\,\pi\,AD\times DB}\\
&\xrightarrow[]{}& {\color{magenta}NOP=AD\times DB}\\
{\color{violet}LG=CD^2} &\xrightarrow[CN2]{}& NOP, LG=AD\times DB, CD^2\\
 {\color{red}NOL, LG=CEFB}  \\
 {\color{blue}CEFB\ is\ CB^2} &\xrightarrow[]{}&  AD\times DB, CD^2= CB^2.
    \end{eqnarray*}

Similarity, we schematize Euclid's proof of the next proposition starting from the conclusion: the rectangle \textit{AM} is equal to the gnomon \textit{NOP}.

II.6 
\pvn\textit{Diorismos}
\[AD\times DB, CB^2=CD^2.\]
\textit{Apodeixis}
\begin{eqnarray*}
... &\xrightarrow[]{}& AM=NOP\\
DM=DB &\xrightarrow[]{}& {\color{violet}AM\,\pi\,AD\times DB}\\
&\xrightarrow[]{}& {\color{magenta}NOP=AD\times DB}\\
{\color{violet}LG=BC^2} &\xrightarrow[]{}& AD\times DB, BC^2=NOP, BC^2\\
 NOL, BC^2=CEFB  \\
 {\color{blue}CEFD\ is\ CD^2} &\xrightarrow[]{}&  AD\times DB, BC^2= CD^2.
    \end{eqnarray*}

In both propositions, figures not represented on the diagrams are to be equal to the squares represented on the diagrams.
 Euclid's job is to show that these not represented are equal to some figures represented on the diagrams. The gnomon \textit{NOP} plays the crucial role in that process.

When one pays no attention to the distinction between figures  in terms of the representation on the diagram, these proofs are alike. However, in II.5, when Euclid takes together the square \textit{LG} and the gnomon \textit{NOP}, they make a figure represented on the diagram. Thus, equality $NOP, LG=CEFB$ is based on  visual evidence. In II.6,  Euclid adds  the square on \textit{BC}, which is not represented on the diagram, to the gnomon \textit{NOP}. Thus,  equality
$NOP, BC^2=CEFB$ can not be based on visual evidence here. In fact, Euclid skips an argument justifying this step.
 Thus, the second proof is  more abstract.

Let us consider the sequence of propositions II.5--8 from the perspective of visible and invisible figures.
In II.5, the equality $NOL, LG=CEFB$ is based on visual evidence. In II.6, the equality $NOL, BC^2=CEFB$ is not so obvious, yet in the company of II.5 it is almost the same. In II.7, Euclid adds to the gnomon \textit{KLM}, the complementing square \textit{DG} and another one  placed on the same diagonal \textit{DB}. These figures are represented on the diagram, yet the equality $KLM, DG, GB=ADEB, GB$ involves overlapping figures. In II.8, Euclid considers overlapping figures but not represented on the diagram.  It is typical of Euclid sequence of micro-steps, similar, e.g. to the first propositions in his theory of equal figures, when he considers parallelograms on the same base, then on equal bases (I.35--36), triangles on the same base, then on equal bases (I.37--38). Therefore, when II.5,\,6 are considered in isolation, they provide almost the same result. When viewed in a bigger picture, they pave a way to a more abstract diagrams. 

%In sum, II.5--6 when viewed in a wider perspective,  play of role in a processes of introducing overlapping figures.
%
%In the next section, we will show that when viewed from the perspective of application in II.11, 14, they differ  more radically. 
 
\subsubsection{Van der Waerden's and Corry's interpretations} Here is van der Waerden's interpretation: ``We see therefore, that, at bottom, II 5 and II 6 are not propositions, but  
solutions of problems; II 5 calls for the construction of two segments \textit{x} and \textit{y} of which the 
sum and product are given, while in II 6 the difference and the product are given. 
The applications in the \textit{Elements} themselves are consistent with this view" (Waerden 1961, 121). 

In his view, II.5 and II.6 are two propositions for one formula, namely $(\frac{x+y}2)^2=xy+(\frac{x-y}2)^2$, 
which is why they need to be interpreted as solutions of problems rather than mere propositions.  To illustrate his idea, van der Waerden interprets proposition II.11 as the solving of a specific equation. Yet, II.11 also allows
 a standard, say a Hilbert-style interpretation, where II.6 is referred to in order to get the result $CF\times FA, AE^2=EF^2$. 
 
 We interpret II.6 as lemma which is applied in II.11, while  II.11 we view as the crucial step in Euclid's construction of dodecahedron -- a regular solid  foreshadowed in Plato's \textit{Timaeus}.

  Corry's interpretation is as follows: ``if we remain close to the Euclidean text we have to admit that, particularly in the cases of II.5 and II.6, both the proposition and its proof are formulated in purely
geometric terms. There are no arithmetic operations involved, and surely there is no
algebraic manipulation of symbols representing the magnitudes involved. The entire
deduction relies on the basic properties of the figures that arise in the initial construction
or that were proved in previous theorems (which in turn were proved in purely
geometric terms)" (Corry 2013, 647).

Our schemes clearly show that Euclid's \textit{deduction} only partly relies on constructed figures.  Euclid's results in Book II also concern  figures which are not represented on diagrams, and not constructed.
\subsection{The use of II.5--6}

We present propositions II.5--6 as lemmas. Indeed, II.5 is applied in II.14, II.6 -- in II.11.
The word pattern of these references is the same:  Euclid simply  repeats the \textit{ekthesis} with new  names of the respective points. In regard to II.5 it is: ``For let any straight-line AB have been cut -- equally at
C, and unequally at D". 
As for II.6 it is: ``For since the straight-line AC has been cut in half at
 E, and FA has been added to it".  In this way, the applied propositions are identified by the patterns of the cut points:
 with II.5 it is equally and unequally, with II.6 -- at half and a line added. This style of references compels Euclid to prove two propositions.  Nevertheless, we provide a detailed analysis of the use of these propositions, as it reveals another relation between visible and invisible figures. In II.1--8, Euclid starts with visible figures to get a relation between invisible figures. In II.11,\,14, by referencing to II.5,\,6, he starts from  invisible figures to get a relation between visible ones. Therefore, when one ignores Euclid's proof techniques, one can still consider propositions II.11,\,14 as a relation between visible figures, and retain a  Euclid drawing of individual lines and circles.
 
% Yet, the difference between II.5 and II.6 is more striking when we consider 
% diagrammatic patterns of their use.
 
 \begin{figure}
\begin{center}
\includegraphics{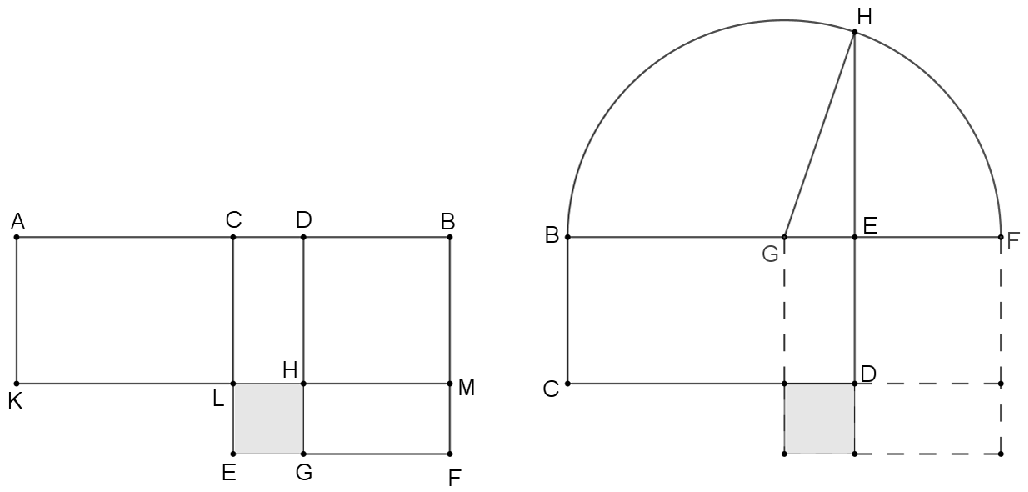}
\caption{The use of II.5 in II.14.} \label{figII5_14b}
\end{center}
\end{figure}

In II.14, it is required to construct a square equal to a rectilinear figure \textit{A}. Due to a triangulation technique, \textit{A} is turned into a rectangle \textit{BCDE}. Here is where our scheme starts.

\begin{eqnarray*}
... &\xrightarrow[II.5]{}& BE\times EF, GE^2=GF^2\\
GF=GH &\xrightarrow[]{}& BE\times EF, GE^2=GH^2\\
HE^2, GE^2=GH^2&\xrightarrow[]{}& BE\times EF, GE^2=HE^2, GE^2 \\
 &\xrightarrow[CN3]{}& BE\times EF=HE^2  \\
   EF=ED   &\xrightarrow[]{}& {\color{violet}BD\,\pi\, BE\times EF}\\
   &\xrightarrow[]{}& {\color{magenta}BD=HE^2}\\
   BD=A&\xrightarrow[]{}& A=EH^2.
    \end{eqnarray*}

When we add the diagram of II.5 to the diagram of II.14,  the proof  supported by the compound diagram goes smoothly (see Fig. \ref{figII5_14b}). On the one hand, by II.5, 
rectangle \textit{BCDE} and the gray square are equal to the square on \textit{GF}, which is equal to the square on \textit{GH}. On the other hand, by I.47, the square on \textit{GH} is equal to the square on \textit{HE} and the gray square. By substitution, $BCED=EH^2$.
The final result concerns figures represented on the diagram, modulo the square on HE represented by its side.

Note, however, that the result was achieved by a reversal of our rule (3).  Specifically, by II.5 and CN3,
$BE\times EF=HE^2$. By substitution rule (1), $BD\,\pi\,BE\times EF$. As a result, $BD=EH^2$. We can turn this process into the following rule
\begin{equation} X=Y\times Z,\ W=Y\times Z \Rightarrow  X=W.\    \end{equation}

The invisible figure, the rectangle contained by \textit{Y}, \textit{Z}, enables the relation between visible figures \textit{X},\textit{W}.

\begin{figure}
\begin{center}
\includegraphics{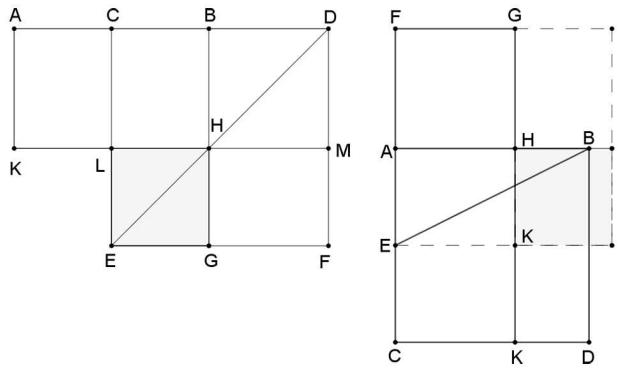}
\caption{The use of II.6 in II.11.} \label{figII6_11}
\end{center}
\end{figure}

Here is the scheme of  II.11 starting from the reference to II.6.

\begin{eqnarray*}
... &\xrightarrow[II.6]{}& CF\times FA, AE^2=EF^2\\
EF=EB &\xrightarrow[]{}& {\color{violet}CF\times FA, AE^2=EB^2}\\
\angle A=\pi/2 &\xrightarrow[I.47]{}& AB^2, AE^2=EB^2 \\
 &\xrightarrow[]{}& CF\times FA, AE^2= AB^2, AE^2 \\
     &\xrightarrow[CN3]{}& CF\times FA= AB^2\\
  AF=FG &\xrightarrow[]{}& {\color{violet}FK\,\pi\,CF\times FA}\\
   {\color{blue}AD\ is\ AB^2}&\xrightarrow[]{}& FK=AD\\
   &\xrightarrow[CN3]{}& FK\setminus AK= AD\setminus AK\\
   &\xrightarrow[]{}& {\color{red}FH=HD}\\
    AB=BD  &\xrightarrow[]{}&  {\color{violet}HD\,\pi\, AB\times BH}\\
   {\color{blue}FH\ is\ HA^2} &\xrightarrow[]{}& {\color{magenta}AB\times BH= HA^2}.
    \end{eqnarray*}

We marked in red the equality $FH=HD$. In fact, equality $FK\setminus AK=FH$, and $AD\setminus AK=HD$ are based on visual evidence. Yet, Euclid skipped these steps.

The rectangle $FCKG$ and the gray square are equal to the square on \textit{EF}, which is equal to the square on \textit{EB}. 
 
Here, by II.6 and CN3, Euclid gets  the equality $CF\times FA= AB^2$. Then, by substitution rules, he aims to turn it into the equality of of figures represented on the diagram,  $FK=AD$. Here,  we can notice the pattern of our rule (4)
\[CF\times FA= AB^2,\ FK\,\pi\, CF\times FA\Rightarrow FK=AB^2.\]

In the rest of the proof, Euclid handles  visible figures, and  the next step relies on visual evidence, $FH=HD$.

Let us have a look at figures Fig. \ref{figII5_14b}, \ref{figII6_11}. In II.14, when we  apply the diagram of II.5 to the line \textit{BF}, no auxiliary lines are needed to finish the proof (modulo the square on HE). In II.11, when  
we adopt the same procedure, the square on \textit{AB}  deforms the diagram  II.6, in a way. 

Finally, let us adopt a mechanical perspective known, for example, through Descartes' drawing instruments; see e.g. (Descartes 1637, 318, 320, 336). Diagram II.11 is, in fact, a project of a machine squaring a rectangle, where a sliding point \textit{E} determines its perimeter. As  figure \textit{A}  can change in the original diagram II.11,  line \textit{GE}  has to change accordingly. In this context, the term \textit{at random}, applied also as a synonym of \textit{unequally}, may suggest a dynamic interpretation. On the other hand, in II.11, the solution is determined by the right-angle triangle \textit{AEB}, and no line can play the role of a variable.

%While from the algebraic point of view there is no difference whether the cut points are at half and at random or at half and a line is added, 

%\begin{figure}
%\begin{center}
%\includegraphics{rysunki/II11}
%\caption{\textit{Elements} II.11.} \label{figII11}
%\end{center}
%\end{figure}

\section{Interpretations of Book II}

 So far, we commented on the recent interpretations of Book II  regarding the specific aspects of our schemes. In this section, we discuss broader analyses.
\subsection{Historians}

\subsubsection{Ken Saito}
 Saito  interprets the propositions of Book II as a relation between visible and invisible figures. He writes:  ``The propositions II 1--10 are those concerning invisible figures,  and they must be proved by reducing invisible figures to visible ones, for one can
apply to the latter the geometric intuition which is fundamental in Greek geometric
arguments" (Saito 2004, 167). 

Instead of the  description  of visible vs invisible, we prefer to address this duality as represented vs not represented on a diagram. It seems a better choice, since two sides of \textit{parallelograms contained by straight-lines containing a right angle} are visible, that is,  represented on a diagram.
Nevertheless,  we should give Saito the  credit for this general observation, especially as 
 Euclid scholars  usually uphold the dogma that ``Greek mathematical proofs are about specific objects in specific diagrams" (Netz 1999, 241).\footnote{In (Błaszczyk, Petiurenko 2020) we identify a tendency in the \textit{Elements} to eliminate visual aspects 
in order to achieve a generality founded on theoretical grounds alone. Thus, II.1 to II.4 exemplify a trend rather than atypical arguments.} Since we identify  the rules relating visible and invisible figures, one may view our study as a  development of Saito's basic observation.

\subsubsection{Leo Corry}
 Corry  applies Saito's distinction of visible vs invisible in his analysis of Book II.  Accordingly, regarding proposition II.1, he formalizes its \textit{diorismos} as the following equation\footnote{Unfortunately, instead of Euclid's \textit{parallelograms contained by},  Corry applies his own term, namely ``R(CD, DH ) means the rectangle
built on CD, DH". It corresponds to our suggestion that Corry pays no attention to the renaming technique characterized above.}
\begin{equation}\tag{Eq. 1} R(A, BC)=R(A,BD)+R(A, DE)+R(A,EC), \end{equation}
and the starting point of Euclid's proof (the formula in red in our scheme) as
\begin{equation}\tag{Eq. 2}  R (BG, BC) = R (BG,BD) + R (DK,DE) + R(EL, EC).      \end{equation} 

Then, he points out that a relation between these equalities  can be explained in terms of  visible and invisible figures.
  Hence, Corry writes: ``what Saito draws our attention to, in particular, is
the fact that the rectangles used in (Eq. 2) are all `visible' in the diagram, whereas
those of (Eq. 1) are `invisible'. [...] In other words, situations
embodied in (Eq. 2) [...] involve visible figures and hence do not require
further justification other than what the figure itself shows. The situation embodied in
(Eq. 1), in contrast, does require a proof precisely because the rectangles involved are,
as indicated by Saito, invisible. In Book II, then, Euclid shows how the properties of
invisible figures can be derived from those of visible ones" (Corry 2013, 650--651).

Regarding the crucial point, namely ``how the properties of
invisible figures can be derived from those of visible ones", Corry's explanation is as follows: 
``The proof itself, on the other hand, is based on (i) taking a segment BG = A, (ii)
constructing the parallelograms and proving on purely geometric grounds (using I.34)
that DK = A = EL, and (iii) then realizing that, according to the diagram:
\begin{equation}\tag{Eq. 2}  R (BG, BC) = R (BG,BD) + R (DK,DE) + R(EL, EC).       \end{equation} 
So, what is the big difference between (Eq. 1) and (Eq. 2) and in what sense does
the latter prove the former? Notice, in the first place, that proving DK = A = EL is
fundamental since otherwise the three rectangles in the figure cannot be concatenated
into a single one in (Eq. 2). But what Saito draws our attention to, in particular, is
the fact that the rectangles used in (Eq. 2) are all `visible' in the diagram, whereas
those of (Eq. 1) are 'invisible'" (Corry 2013, 650).

Indeed, step (i) is the \textit{kataskeu\={e}} part of Euclid's proof.   As for  step (ii), it is Corry's argument rather than Euclid's, since the  text of the  proposition is: ``DK, that is to say BG, is equal A". It means that Euclid does not justify the equalities $DK=BG=A$. Nevertheless, it is a favorable argument, if needed.   Step (iii) is what we consider as  visual evidence. However, steps (i)--(iii) do not provide a complete account of Euclid's proof.    Corry does not explain how Euclid relates the visible figure $R(BG, BC)$ and the invisible  $R(A, BC)$. 
The simple observation that,  on the one hand, there are  visible figures, on the other hand,  invisible ones, does not tell us how Euclid turns  equation Eq.1 into  equation Eq.2. We believe that our substitution rules enable 
an adequate explanation.

\subsubsection{Ian Mueller}  In most of his review of Book II, Mueller argues against algebraic interpretation; see  (Mueller 2006, 41--52, 301--302). In this subsection, we try to separate  his own interpretation from this polemic.
 In section \S\,3, we have  shown that Mueller adopts a notation  which revokes the distinction between visible and invisible figures. Let us recall his definitions:
``\textbf{O}(x,y) is used to designate  a rectangle with arbitrary straight lines equal to $x$ and $y$  as adjacent sides", ``I use \textbf{T}(\textit{x}) to stand for the square on a straight line equal to $x$" (Mueller 2006, 42, 45) 

 Actually, there is no significant difference between $\textbf{O}(x,y)$ and the  algebraic term $xy$.
On the one hand, algebraic interpretation takes it for granted that $x(y+z)=xy+xz$, on the other hand, the rule 
$\textbf{O}(x,y+z)\simeq\textbf{O}(x,y)+\textbf{O}(x,z)$ 
 is   self-evident for Mueller.  Moreover, like in algebraic interpretation, the term $\textbf{O}(x,y)$ is applied to visible and invisible figures in the same way.
 
For example, here is Mueller's reading of II.1: ``It should be clear, once the construction is described, II.1 becomes a geometrically trivial proposition" (Mueller 2006, 42). However, according to our scheme of II.1, only the first step represented by the formula in red is trivial. The rest of the proof is far from obvious. 

As long as Mueller interprets the \textit{diorismos} parts, his formalism  works well. When he seeks to analyze Euclid's proofs,  it leads him astray. Here is his reading of the proof of II.4:\footnote{See (Mueller 2006, 45--48).} 

$\textbf{T}(BC)\simeq   square  BG$,

$\textbf{T}(AC)\simeq   square \ HF$,

$\textbf{O}(AC,BC)\simeq  rectangle\ AG\simeq [by\ I.43]\ rectangle \ GE$,

\textit{since} $\textbf{T}(AB)\simeq   square \ AE$, \textit{the theorem follows, that is }

$\textbf{T}(AB)\simeq  \textbf{T}(AC)+ \textbf{T}(CB)+ 2 \textbf{O}(AC, CB)$.

\pvn

The relations between visible and invisible figures, which we explain via substitution rules,  are covered by the congruence $\simeq$ alone in Mueller's interpretation. 
Thus, the line of arguments 
$$\textbf{O}(AC,BC)\simeq  rectangle\ AG\simeq [by\ I.43]\ rectangle \ GE,$$ 
aims to interpret Euclid's two different relations:  AG \textit{is contained by} AC, BC, and   $AG=GE$.
 Moreover, $\simeq$ is used to explain what we call the renaming of figures. Thus,  Euclid's ``CGKB is the square on CB", 
  Mueller interprets also by the  congruence: $\textbf{T}(BC)\simeq   square  BG$. Since this congruence is supposed to be transitive -- Mueller does not explain why  it is so, in the context of Book II -- Euclid's proof seems to go smoothly. However, it breaks as the conclusion $\textbf{T}(AB)\simeq  \textbf{T}(AC)+ \textbf{T}(CB)+ 2 \textbf{O}(AC, CB)$ comes out of nothing. Mueller skips any reference to the relation $AE\simeq \textbf{T}(AC)+ \textbf{T}(CB)+ 2 \textbf{O}(AC, CB)$. Why?
  
  Euclid's argument is this:

\begin{eqnarray*}
 {\color{red}HF,\, CK,\, AG,\, CE=ADEB}, &&\\
 {\color{blue}ADEB\ is\ AB^2} &\xrightarrow[]{}& {\color{magenta}AB^2=AC^2,\, CB^2,\, 2 AC\times CB}.
    \end{eqnarray*}

  Since the term $\textbf{O}(x,y)$ applies both to visible and invisible figures, within Mueller's formalization, it would
  assume such a form:
  
  Since
  \begin{eqnarray*}
   \textbf{T}(AC)+ \textbf{T}(BC)+\textbf{O}(AC, CB)+\textbf{O}(AC, CB)\simeq square\ AE, &&\\
 {\textbf{T}(AB)\simeq square\ AE},\end{eqnarray*}
 
 then
 \begin{eqnarray*}
 \textbf{T}(AB)\simeq \textbf{T}(AC)+ \textbf{T}(CB)+ 2 \textbf{O}(AC, CB).
   \end{eqnarray*}

It would result in a vicious circle argument. Therefore, Mueller had to skip Euclid's reference to visual evidence. As a result, he mischaracterized Euclid's proof.

\subsection{Mathematicians}

\subsubsection{Bart van der Waerden}

Van der Waerden is a prominent advocate for the so-called geometric algebra interpretation of Book II.  Recent papers  by Victor Bl\aa sj\"{o} and Mikhail Katz  recount this fascinating debate between mathematicians and historians.\footnote{See (Bl\aa sj\"{o} 2016), (Katz 2020)} From our perspective, however, it is too abstract, as it does not stick to source texts  closely enough.  The analysis of van der Waerden's arguments below certyfies our claim.  

Van der Waerden writes: ``When one opens Book II of the \textit{Elements}, one finds a sequence of propositions which are nothing but  geometric formulations of algebraic  rules. So, e.g., II.1:  [...] corresponds to the formula $a(b+c+...)=ab+ac+...$. II 2 and 3 are special cases of this proposition. II 4 corresponds to the formula 
$(a+b)^2=a^2+b^2+2ab$. The proof can be read off immediately from Fig 34" (Waerden 1961,  118). 

\begin{figure}
\begin{center}
\includegraphics{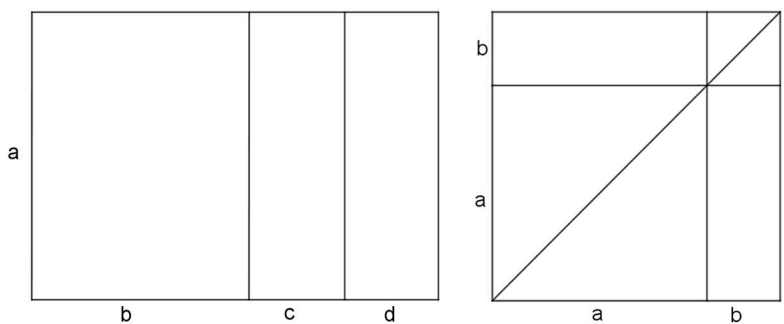}
\caption{Van der Waerden's diagrams for  II.1 and II.4.} \label{figVdW}
\end{center}
\end{figure}

Our Fig. \ref{figVdW} represents van der Waerden's Figures 33 and 34. They aim to  emulate Euclid's diagrams accompanying  propositions II.1  and II.4. 
Let us notice that these diagrams differ from Euclid's in regard to the names of line segments -- it never happens in the \textit{Elements} diagrams that different individual lines have got the same name. Moreover, there is no counterpart of line A in  Figure 33.  It looks like van der Waerden had to modify Euclid's diagrams to develop his interpretation.

Now, let us take the formula  $a(b+c+d)=ab+ac+ad$ designed to correspond  to proposition II.1. Which part of the proposition does it formalize: the \textit{diorisomos}, or the starting point of the proof (the formula in red in our scheme)? In fact, since van der Waerden's diagram does not represent  line A,  his account of Euclid's proof  would  look like this
\begin{eqnarray*}
a(b+c+d)=ab+ac+ad &\xrightarrow[]{}& \\
 &\vdots & \\
  &\xrightarrow[]{}& a(b+c+d)=ab+ac+ad.
    \end{eqnarray*}
    
Thus, there is no need for any proof at all. 

The same applies to his interpretation of proposition II.4. 
Instead of 
\begin{eqnarray*}
 {\color{red}HF,\, CK,\, AG,\, CE=ADEB} &\xrightarrow[]{}&\\
&\vdots&\\
  &\xrightarrow[]{}& {\color{magenta}AB^2=AC^2,\, CB^2,\, 2 AC\times CB},
    \end{eqnarray*}

van der Waerden-style proof would look like this
\begin{eqnarray*}
 (a+b)^2=a^2+b^2+2ab &\xrightarrow[]{}&\\
&\vdots&\\
  &\xrightarrow[]{}& (a+b)^2=a^2+b^2+2ab.
    \end{eqnarray*}

There is also no need for any proof. Indeed, regarding proposition II.4, he writes ``The proof can be read off immediately from Fig 34".
However, in the \textit{Elements}, only the equality in red is read off the diagram, while the final conclusion requires some arguments.

In sum, whatever van der Waerden interprets, these are not Euclid's propositions.

Finally,  we find the following  speculations:  ``We were not able to find any interesting geometrical problem that would give rise to theorems like II 1--4. On the other hand, we found that the explanation of these theorems as arising from algebra worked well" (Waerden 1975, 203). There is no  need  to dispute  whether the distinction between figures represented and not represented on a diagram is a geometrical problem. Whatever it is, it provides an explanation for Euclid's propositions II.1 to II.4. As regards strictly \textit{geometrical problems}, II.4 is applied in II.12 which is the ancient counterpart for the cosine rule for obtuse triangle. Below diagram illustrates the use of II.4

\begin{figure}
\begin{center}
\includegraphics{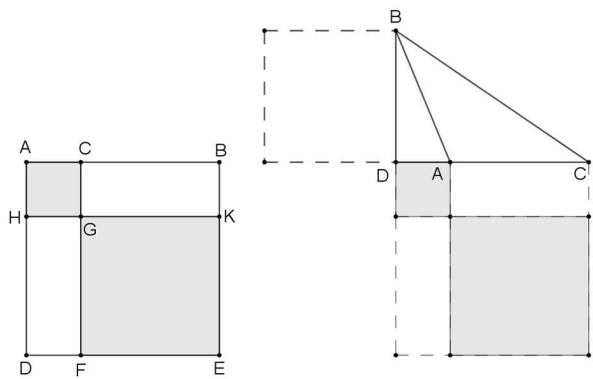}
\caption{Euclid, \textit{Elements}, II.4, 12.} \label{figII4_12}
\end{center}
\end{figure}

Here is the  respective scheme.

\pvn\textit{Diorismos}
\[CB^2=CA^2, AB^2, 2CA\times AD.\]
\textit{Apodeixis}
\begin{eqnarray*}
 &\xrightarrow[II.4]{}&  DC^2=AC^2,AD^2, 2CA\times AD\\
 &\xrightarrow[]{}& DC^2, DB^2=AC^2,AD^2, DB^2, 2CA\times AD\\
 \angle D=\pi/2 &\xrightarrow[]{}& CD^2, DB^2=CB^2\\
 &\xrightarrow[]{}&AD^2, DB^2=AB^2\\
 &\xrightarrow[]{}& CB^2=CA^2, AB^2, 2CA\times AD.
    \end{eqnarray*}

The reference to II.4 is easily identified by the way the line CD is cut at the point A: ``since the straight-line CD has been cut, at random, at point A, the one on DC is thus equal to  the  squares on CA, AD, and twice the
 rectangle contained by CA, AD".\footnote{The $\sin$ of the angle DAB, i.e., $-\cos$ of the angle BAC, is the fraction $AD/BA$. That is how we get the modern version of the cosine rule.}
Again, we will not dispute whether the cosine rule is an ``interesting geometrical problem".

In section \S\,6.1.1, we showed that van der Waerden interprets 
II.4,\,5,\,11 as solving specific equations. Significantly, he did not provided similar interpretations for II.12--14.

\subsubsection{John T. Baldwin and Andreas Mueller}

 (Baldwin, Mueller 2019) provides a series of arguments for autonomy of geometry. It includes historical, conceptual and model theoretical ones. As regards history, the paper develops a geometric interpretation of Book II as opposed to van der Waerden's `geometric algebraic' interpretation, as they call it. Geometry in this context, implicitly, means for them a  study of figures represented on the diagrams.  

Baldwin and Mueller place Book II within Euclid's theory  of equal figures:
``On reflection, there is a natural geometric motivation for the main themes of
Book II: Determine a precise method for determining which of two disjoint rectilinear
figures (polygons) has the greater area" (Baldwin, Mueller 2019, 8). In fact, only two propositions of Book II, namely II.5 and II.14, complete the theory of equal figures as developed in Book I.

Baldwin and Mueller continue: ``Thus, Proposition II.2 certainly implies that if a square is split into two non
overlapping rectangles the sum of the areas of the rectangles is the area of the square" (Baldwin, Mueller 2019, 8). However, what they refer to it is the starting point of II.2, not the conclusion. This starting point (the formula in red, in our scheme), as based on visual evidence needs no proof. 

Accordingly,  they present II.5 as a dissection proof; see (Baldwin, Mueller 2019, 9--10).
Here is their proof schematized according to the rules we have already applied to Euclid's proofs.

By construction, $DB\cong BM$, and $CD\cong MF$. The rest is as follows.
\begin{eqnarray*}
{\color{red}ADHK=ACLK, CDHL} &\xrightarrow[]{}& \\
ACLK\cong BFGD &\xrightarrow[]{}& \\
{\color{blue}CB^2\ (is)\ (CBFE)} &\xrightarrow[]{}& \\
 {\color{red}CB^2=BFGD, CDHL, LHGE}\\
   LHGE=CD^2 &\xrightarrow[]{}&.
    \end{eqnarray*}

Indeed, Baldwin--Mueller's proof is a series of observations rather than arguments. It also does not provide a final conclusion.  To get it, one should apply
the substitution rule, namely ${\color{magenta}CB^2=ADHK, CD^2}$. 

In the above scheme, equalities in red interpret the phrase ``is composed", the formula $LHGE=CD$ interprets the phrase ``LHGE (which has the same area as the square on CD)". Thus, Baldwin and Mueller provide a  styling on Euclidean proof  rather than an interpretation of the actual Euclid's proof.

Historians often point out that algebraic interpretation ignores the role of gnomons in Book II. Baldwin and Mueller
managed to turn that objection into a more specific argument, namely:
``Much of Book II considers the relation of the areas of various rectangles,
squares, and gnomons, depending where one cuts a line. While gnomons have
a clear role in decomposing parallelograms, the algebraic representation for the
area of gnomon, is not a tool in polynomial algebra. That is, while such equations
as $(a+b)(a-b)=a^2-b^2$ or the formula for product of binomials are tools
in algebra which have nice geometric explanation, the area of a gnomon has an
algebraic expression, $2ab+a^2$, which does not recur in algebra (e.g., as a method
of factorization)" (Baldwin, Mueller 2019, 9).

Although Baldwin and Mueller emphasize the role of gnomons, in fact, in their proof of II.5, Euclid's gnomon \textit{NOP} is simply a composition of two rectangles: \textit{BFGD}, \textit{CDHL}.  As a result, they do not provide a counterpart of Euclid's decisive argument, namely $NOP=AD\times DB$.  

What is, then, the role of the gnomon in II.5.
Starting from the equality $NOP=AD\times DB$, with \foreignlanguage{polutonikogreek}{proske'ijw}, Euclid  refers to CN2 to get $NOP, LG=AD\times DB, CD^2$.
What Baldwin and Mueller get by visual evidence, Euclid gets by deduction.

How about Baldwin--Mueller congruence $ACLK\cong BFGD$? By no means it is obvious, as we gave up an algebraic mode.  Dissection also seems useless. 

Here is how Euclid gets the result $AH=NOP$ in proposition II.5 (see Fig. \ref{figII5}).

\begin{eqnarray*}
&\xrightarrow[I.43]{}&CH=HF\\
&\xrightarrow[CN2]{}& CH, DM=HF,DM\\
&\xrightarrow[]{}& {\color{red}CM=DF}\\
AC=CB&\xrightarrow[]{}&CM=AL\\
AL,CH=CM,HF&\xrightarrow[]{}& {\color{red}AH=NOP}.
\end{eqnarray*}

While Baldwin and Mueller did not manage to represent Euclid's reliance on gnomons in II.5, contrary to Euclid, they apply gnomon in their proof of II.14. 

In II.14, let us remind,  it is required to construct a square equal to a rectilinear figure $A$. Within the theory of equal figures,  $A$ is turned into a rectangle $BCDE$. And here is where our scheme starts.

\begin{eqnarray*}
... &\xrightarrow[II.5]{}& BE\times EF, GE^2=GF^2\\
GF=GH &\xrightarrow[]{}& BE\times EF, GE^2=GH^2\\
HE^2, GE^2=GH^2&\xrightarrow[]{}& BE\times EF, GE^2=HE^2, GE^2 \\
 &\xrightarrow[CN3]{}& BE\times EF=HE^2  \\
   EF=ED   &\xrightarrow[]{}& {\color{violet}BD\,\pi\, BE\times EF}\\
   &\xrightarrow[]{}& {\color{magenta}BD=HE^2}\\
   BD=A&\xrightarrow[]{}& A=EH^2.
    \end{eqnarray*}
\begin{figure}
\begin{center}
\includegraphics{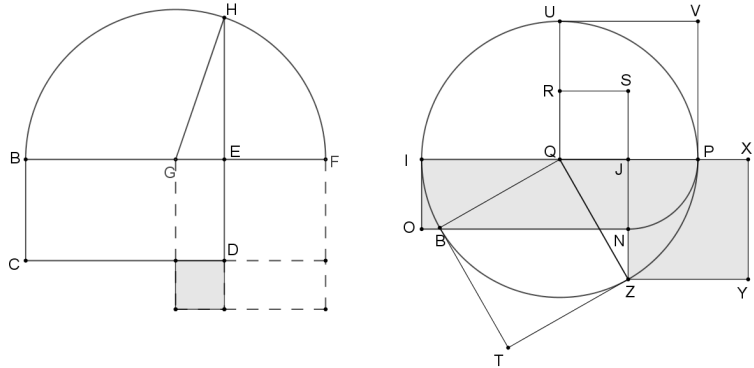}
\caption{Euclid's proof of II.14 (left), Baldwin-Meuller version (right).} \label{figBM}
\end{center}
\end{figure}

And here is Baldwin--Meuller proof (see Fig. \ref{figBM}).

\begin{eqnarray*}
... &\xrightarrow[II.5]{}& IJNO, QJ^2=QP^2\\
QP\cong QZ &\xrightarrow[]{}& {\color{magenta}QP^2=QZ^2}\\
&\xrightarrow[I.47]{}& QP^2=QJ^2, JZ^2 \\
 &\xrightarrow[CN3]{}& QP^2\setminus QJ^2=UVPJR=IJNO=JZ^2.
    \end{eqnarray*}

The reference to the gnomon $UVPJR$ does not add anything to this proof -- from the deductive perspective  it  is useless. Yet, Baldwin and Mueller created a diagram for II.14 in which every argument (every line in the scheme of their proof) is represented by an individual figure. 
Yet, they did not manage to find a similar solution for this argument
\[QP\cong QZ \rightarrow QP^2=QZ^2.\]

  Now, let us compare Baldwin--Mueller's and Euclid's diagrams. On the one hand, there is a complex composition of rectangles and squares designed to represent  every textual argument (Baldwin and Mueller's arguments, instead of Euclid's). On the other hand, there is a simple composition  which represents, actually, abstract arguments.  If we take Euclid's diagram in its original form (see Fig. \ref{figII14}), we can still recognize how II.5 is applied, and the diagram becomes an ideogram.  If there are no other  reasons for the use of the terms \textit{rectangle contained by}  and \textit{square on},   
  this diagrammatic representation of  Euclid's abstract arguments is the most solid one.
   
%With Baldwin, Mueller  Eculid's proof should be schematized like this
%\begin{eqnarray*}
%... &\xrightarrow[II.5]{}& IJNO, QJ^2=QP^2\\
%&\xrightarrow[]{}& IJNO, QJ^2 =QJ^2, JZ^2 \\
% &\xrightarrow[]{}& IJNO=JZ^2.
%    \end{eqnarray*}
%This simple proof demonstrates how Euclid applies I.47.

%\subsection{Metaphysical vs pragmatic interpretation}

\section{Descartes' lettered diagrams}

In this final section, we discuss one of Descartes' diagrams included in his \textit{Geometry}. By this contrast, we aim to  
underline  the crucial steps in  Euclid's procedures used in Book II.

In the   fragment cited below, Descartes illustrates one of his rules for solving   fourth degree equations. The details of his problem are irrelevant for our analysis, therefore they are skipped.
 \begin{figure}
\begin{center}
\includegraphics{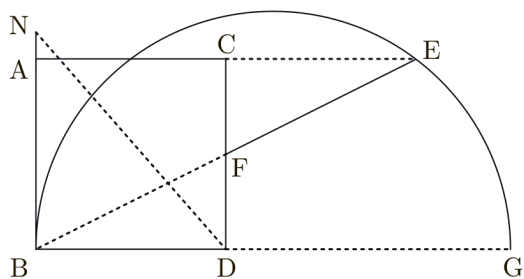}
\caption{Descartes, \textit{Geometry}, p. 388.} \label{figK}
\end{center}
\end{figure}

Here is Descartes' account of the diagram represented in Fig. \ref{figK}:
\pvn ``For, putting $a$ for BD or CD, and $c$ for EF, and $x$ for DF, we have
CF $=a-x$, and, as CF or $a-x$ is to FE  or $c$, so FD or $x$ is to BF,
which as a result gives $\frac{cx}{a-x}$. Now, in the right-angle triangle BDF one side is $x$,
another is $a$, then their squares, which are (\textit{leures carr\'es, qui sont}) $xx+aa$,  are equal to the square of base, which is $\frac{ccxx}{xx-2ax+aa}$"  (Descartes 1637/2007,  388/191).\footnote{Translation after (Descartes 2007) slightly modified.}

On this diagram, Descartes applies the ancient technique of naming vertices and intersections of lines with capital letters. Nevertheless, when it comes to analysis, he assigns  new letters, specifically one letter to different lines: in this case, he puts $a$ for \textit{BD} and \textit{CD}. We have already observed that Euclid applies the same trick in proposition II.1. Yet, in the following propositions, he does not explore this idea. As for Descartes, it is his standard procedure.\footnote{Another spectacular example is Descartes' analysis of Euclid's proposition III.36; see (Descartes 1637, 302). While showing how to construct roots of the second degree polynomial, he assigns  the name $\frac{1}{2}a$ to lines ON and NL. Yet, that assignment is implicit.} 
 We believe that is where the algebra begins:  giving the same name to different objects. 
  The next step in algebraic account of geometry is symbolic computation. As we all know, Descartes developed a technique of symbolic computation. It is the arithmetic of line segments. However, in regard to this      
initial step,  it seems that Descartes does not appreciate it.  
That is how he explains his technique of naming lines:

\pvn ``Finally, so that we may be sure to remember the names of these lines, 
a separate list should always be made as often as names are assigned 
or changed. For example, we may write, AB=1, that is AB is equal 
to 1; GH=a; BD = b, and so on. 
If, then, we wish to solve any problem, we first suppose the solution 
already effected, and give names to all the lines that seem needful for 
its construction, to those that are unknown as well as to those that 
are known.  Then, making no distinction between known and unknown 
lines, we must unravel the difficulty in any way that shows most naturally the relations between these lines, until we find it possible to express 
a single quantity in two ways. This will constitute an equation, since 
the terms of one of these two expressions all together equal to the 
terms of the other" (Descartes 2007, 6--9).

The same letter, \textit{x}, \textit{y} or \textit{z}, can name   \textit{unknown}, and unequal lines.  Yet, Descartes does not mention that the same letter can name different, though known and equal lines. It is his practice, not an explicit method.
From his perspective, the most important notion is that ``a single quantity", i.e. a line segment, can get two names. In his system, these two names make an  equation of two arithmetic terms. However, within the theory of polynomials developed in Book III, these  names can be transformed into a polynomial equation $f(x)=0$, and their original meaning, i.e. a reference to a  specific line, evaporates.

From our perspective, the crux of Descartes' method  consists in  giving one name to different objects: that is the starting point of his equations. 
However, Euclid has a technique of asserting different names to the same object (renaming, in our terms). The transformations of these names allowed him to relate figures which are represented and not represented on diagrams.
The latter  only have one name: it can be  \textit{rectangle contained by} or \textit{square on}. Thus, if one looks for an algebraic prelude in Euclid's Book II, it could be  proposition II.1.
Yet, the technique of assigning one name to different objects has not been developed further in the \textit{Elements}.

%To see the difference, Descartes figure is concrete while its description -- abstrat. Euclid's figure is abstract while description concrete.
%
%In our analysis, in is a relation between what is represented on a diagram a imagined figures. Imagined figures are not represented although their sides are.

\end{document}